\documentstyle[twoside,12pt]{article}
\newcommand\text[1]{{\mbox{ #1 }}}
\newif\ifpctex
\pctexfalse
\ifpctex\hoffset=-lin\voffset=-lin\fi
\ifpctex\def\bb{bf}\else\newfont{\bb}{msbm10 scaled\magstep1}\fi

\newcommand{\A}{{\cal A}}

\newcommand{\F}{{\cal F}}

\renewcommand{\H}{{\cal H}}
\newcommand{\I}{{\cal I}}

\renewcommand{\O}{{\cal O}}

\newcommand{\ph}{\varphi}

\def\nz{\mbox{\bb N}}
\def\rz{\mbox{\bb R}}

\def\cz{\mbox{\bb C}}

\newcommand{\limn}{\lim_{n\to\infty}}

\newcommand{\liminfn}{\liminf_{n\to\infty}}

\newcommand{\qed}{\hfill$\Box$\linebreak\medskip\par}

\newcommand{\pr}{{\em Proof:\quad}}

\newcommand{\be}{\begin{enumerate}}
\newcommand{\ee}{\end{enumerate}}
\newcommand{\bdm}{\begin{displaymath}}
\newcommand{\edm}{\end{displaymath}}
\newcommand{\beq}{\begin{equation}}
\newcommand{\eeq}{\end{equation}}
\newcommand{\beqa}{\begin{eqnarray}}
\newcommand{\eeqa}{\end{eqnarray}}
\newcommand{\beqas}{\begin{eqnarray*}}
\newcommand{\eeqas}{\end{eqnarray*}}

\newcommand{\alternative}[4]
{\left\{
\begin{array}{llll}
#1&#2\\#3&#4
\end{array}
\right. }

\newtheorem {lemma}{Lemma}[section]
\newtheorem {theo}{Theorem}[section]
\newtheorem {bemerkung}{Remark}[section]
\newtheorem {remark}{Remark}[section]
\newtheorem {prop}{Proposition}[section]

\newtheorem {beispiel} {Example}[section]

\newcommand{\sect}[1]{\section{#1}\setcounter{equation}{0}}
\renewcommand{\rho}{\varrho}
\renewcommand{\c}[1]{c_{S_{#1}}}

\newcommand{\f}[1]{f^{S_{#1}}}

\newcommand{\el}{\frac 1{\ell}}

\newcommand\proof{\pr}

\newcommand\sumai[1]{{\sum_{i=#1}^\infty a_i}}
\newcommand\sumaI[1]{{\sum_{i=0}^{#1} a_i}}
\newcommand\nought{\scriptscriptstyle0}
\newcommand\rhop{\rho_{\scriptscriptstyle1}}
\newcommand\rhom{\rho_{\nought}}
\newcommand\Kp{K_{\scriptscriptstyle1}}

\newcommand\cSn[1]{c_{S_{#1}}}
\newcommand\cmSn[1]{z_n}
\newcommand\fSn[1]{f^{S_{#1}}}
\newcommand\fmSn[1]{f^{S_{#1}-1}}

\newcommand\dis[2]{{\bf{#1}}{#2}}
\newcommand\st{\,\,;\,\,}
\newcommand\intersect{3}

\begin{document}
\title
{Absorbing Cantor sets in dynamical systems: Fibonacci maps}
\author{
Henk Bruin, Univ. of Delft, the Netherlands\\
Gerhard Keller, Univ. of Erlangen, FRG\\
Tomasz Nowicki, Univ. of Warsaw, Poland
\thanks{supported by the DFG, NWO, KBN-GR91}
\\
Sebastian van Strien, Univ. of Amsterdam, the Netherlands
\thanks{Part of this work was done at Stony Brook. SvS would like to
thank M. Lyubich and F. Tangerman for some useful discussions.}
}
\date{\relax}
\maketitle
\thispagestyle{empty}
\def\IMSmarkvadjust{0 pt}
\def\IMSmarkhadjust{0 pt}
\def\IMSmarkhpadding{0 pt}
\def\IMSpubltext{Published in modified form:}
\def\SBIMSMark#1#2#3{
 \font\SBF=cmss10 at 10 true pt
 \font\SBI=cmssi10 at 10 true pt
 \setbox0=\hbox{\SBF \hbox to \IMSmarkhpadding{\relax}
                Stony Brook IMS Preprint \##1}
 \setbox2=\hbox to \wd0{\hfil \SBI #2}
 \setbox4=\hbox to \wd0{\hfil \SBI #3}
 \setbox6=\hbox to \wd0{\hss
             \vbox{\hsize=\wd0 \parskip=0pt \baselineskip=10 true pt
                   \copy0 \break%
                   \copy2 \break%
                   \copy4 \break}}
 \dimen0=\ht6   \advance\dimen0 by \vsize \advance\dimen0 by 8 true pt
                \advance\dimen0 by -\pagetotal
	        \advance\dimen0 by \IMSmarkvadjust
 \dimen2=\hsize \advance\dimen2 by .25 true in
	        \advance\dimen2 by \IMSmarkhadjust

%
%
  \openin2=publishd.tex
  \ifeof2\setbox0=\hbox to 0pt{}
  \else 
     \setbox0=\hbox to 3.1 true in{
                \vbox to \ht6{\hsize=3 true in \parskip=0pt  \noindent  
                {\SBI \IMSpubltext}\hfil\break
                ``Wild Cantor Attractors Exist'', {\it Annals of Math.} {\bf 43} (1996),
pp. 97--130.
 
                \vfill}}
  \fi
  \closein2
  \ht0=0pt \dp0=0pt
 \ht6=0pt \dp6=0pt
 \setbox8=\vbox to \dimen0{\vfill \hbox to \dimen2{\copy0 \hss \copy6}}
 \ht8=0pt \dp8=0pt \wd8=0pt
 \copy8
 \message{*** Stony Brook IMS Preprint #1, #2. #3 ***}
}

\SBIMSMark{1994/2}{March 1994}{}
%

%
\ifx\beginpic\undefined\else \fi

\chardef\oldatcat=\the\catcode`\@
\catcode`\@=11

\newskip\hsssglue \hsssglue=0pt plus 1fill minus 1fill \def\hsss{\hskip\hsssglue}

\newdimen\unitlength \newdimen\linethickness
\newdimen\@picheight \newdimen\@xdim \newdimen\@ydim \newdimen\@len \newdimen\@save
\newcount\@multicount \newcount\@xarg \newcount\@yarg
\newbox\@picbox \newbox\@mpbox

\font\tenln=line10     \font\tenlnw=linew10
\font\tencirc=lcircle10 \font\tencircw=lcirclew10
\font\smallln=linew10 scaled 483 

\def\thinlines{\let\linefont=\tenln \let\circlefont=\tencirc
  \linethickness=\fontdimen8\linefont}
\def\thicklines{\let\linefont=\tenlnw \let\circlefont=\tencircw
  \linethickness=\fontdimen8\linefont}
\thinlines

\def\beginpic(#1,#2)(#3,#4){\@picheight=#2\unitlength
  \setbox\@picbox=\hbox to#1\unitlength\bgroup\let\line=\@line
    \kern-#3\unitlength \lower#4\unitlength\hbox\bgroup\ignorespaces}
\def\endpic{\egroup\hss\egroup
  \ht\@picbox=\@picheight \dp\@picbox=\z@
  \leavevmode\box\@picbox}

\def\put(#1,#2)#3{\raise#2\unitlength\rlap{\kern#1\unitlength #3}\ignorespaces}

\def\multiput(#1,#2)(#3,#4)#5#6{\@multicount=#5
 \@xdim=#1\unitlength \@ydim=#2\unitlength \setbox\@mpbox=\hbox{#6}%
 \loop\ifnum\@multicount>0
   \raise\@ydim\rlap{\kern\@xdim \unhcopy\@mpbox}%
   \advance\@xdim#3\unitlength \advance\@ydim#4\unitlength
   \advance\@multicount\m@ne \repeat\ignorespaces}

\def\makebox(#1,#2)#3{\setbox\@picbox=\hbox to#1\unitlength{\hss#3\hss}%
  \@ydim=\ht\@picbox \advance\@ydim-\dp\@picbox
  \ht\@picbox=#2\unitlength \dp\@picbox=\z@
  \leavevmode\lower.5\@ydim\box\@picbox}

\newif\ifneg
\def\@line(#1,#2)#3{\@xarg=#1 \@yarg=#2 \@len=#3\unitlength \leavevmode
 \ifnum\@xarg<0 \reverseline \else \negfalse \@ydim=\z@\fi
 \ifnum\@xarg=0 \@vline
 \else\ifnum\@yarg=0 \@hline \else\@sline\fi\fi
 \ifneg\kern-\@len\else\@save=\@ydim\fi}
\def\reverseline{\negtrue \kern-\@len \@xarg=-\@xarg
 \@ydim=\@len \multiply\@ydim\@yarg \divide\@ydim\@xarg \@yarg=-\@yarg}

\def\@hline{\vrule height.5\linethickness depth.5\linethickness width\@len}
\def\@vline{\kern-.5\linethickness\vrule width\linethickness
  \ifnum\@yarg<0 height\z@ depth\else depth\z@ height\fi\@len
  \kern-.5\linethickness}

\def\@sline{\setbox\@picbox=\hbox{\linefont \count@=\@xarg \multiply\count@ 8
 \ifnum\@yarg>0 \advance\count@\@yarg \advance\count@-9
 \else \advance\count@-\@yarg \advance\count@ 55 \fi \char\count@}%
 \ifnum\@yarg<0 \@picheight=-\ht\@picbox \advance\@ydim\@picheight
 \else \@picheight=\ht\@picbox \fi
 \@xdim=\wd\@picbox \@save=\@ydim
 \loop\ifdim\@xdim<\@len \raise\@ydim\copy\@picbox
  \advance\@xdim\wd\@picbox \advance\@ydim\@picheight \repeat
 \advance\@xdim-\@len \kern-\@xdim
 \multiply\@xdim\@yarg \divide\@xdim\@xarg \advance\@ydim-\@xdim
 \raise\@ydim\box\@picbox}

\def\vector(#1,#2)#3{\@line(#1,#2){#3}%
 \ifnum\@xarg=0 \@vvector \else\ifnum\@yarg=0 \@hvector \else\@svector\fi\fi}
\def\@hvector{\ifneg\rlap{\linefont\char27}\else
 \smash{\llap{\linefont\char45}}\fi} 
\def\@vvector{\ifnum\@yarg<0 \raise-\@len\rlap{\linefont\char63}%
 \else\setbox\@picbox=\rlap{\linefont\char54}\advance\@len-\ht\@picbox
 \raise\@len\box\@picbox\fi}

\def\@svector{\setbox\@picbox=\hbox to\z@{\linefont
 \ifnum\@yarg<0 \count@=55 \@yarg=-\@yarg \else\count@=-9 \fi
 \ifneg\multiply\@xarg16 \multiply\@yarg2
 \else\hss 
  \ifnum\@xarg>2 \multiply\@xarg9 \multiply\@yarg2 \advance\count@29
  \else\ifnum\@yarg>2 \multiply\@xarg16 \multiply\@yarg9 \advance\count@-20
   \else\multiply\@xarg24 \multiply\@yarg3 \fi\fi\fi
  \advance\count@\@xarg \advance\count@\@yarg \char\count@
  \ifneg\hss\fi}
 \raise\@save\box\@picbox}

\def\disk#1{\@len=#1\unitlength \count@='160 \@diskcirc}
\def\circle#1{\@len=#1\unitlength \count@='140 \@diskcirc}
\def\@diskcirc{\setbox\@picbox=\hbox{\circlefont\char\count@}\@xdim=\wd\@picbox
 \leavevmode \ifdim\@len>15.499\@xdim \@bigdc \else \@smalldc\fi}
\def\@bigdc{\ifnum\count@<'160 \@bigcirc
 \else \@len=15\@xdim \@diskcirc\fi}
\def\@smalldc{{\advance\@len-.5\@xdim
 \loop\ifdim\@xdim<\@len \advance\count@\@ne \advance\@xdim\wd\@picbox\repeat
 \hbox{\circlefont\char\count@}}}
\def\@bigcirc{{\circlefont\count@=15
 \setbox\@picbox=\hbox{\char\count@}\@xdim=\wd\@picbox
 \ifdim\@len>2.5\@xdim \@len=2.5\@xdim\fi
 \advance\@len-.125\wd\@picbox
 \loop\ifdim\@xdim<\@len \advance\count@ 4 \advance\@xdim.25\wd\@picbox\repeat
 \@ydim=.5\@xdim \advance\@ydim.5\linethickness
 \setbox\@picbox=\vbox{\hbox{\char\count@\advance\count@-3\char\count@}%
  \nointerlineskip
  \hbox{\advance\count@\m@ne\char\count@\advance\count@\m@ne\char\count@}}%
 \kern-\@ydim\lower\@ydim\box\@picbox}}

\newif\ifovaltl \newif\ifovaltr \newif\ifovalbl \newif\ifovalbr
\ovaltltrue \ovaltrtrue \ovalbltrue \ovalbrtrue
\def\oval(#1,#2){\@xdim=#1\unitlength \@ydim=#2\unitlength
 {\circlefont \setbox\@picbox=\hbox{\char0}
 \ifdim\@xdim<\wd\@picbox \@xdim=\wd\@picbox\fi
 \ifdim\@ydim<\wd\@picbox \@ydim=\wd\@picbox\fi
 \@save=\@xdim \ifdim\@ydim<\@save \@save=\@ydim \fi
 \count@=39
 \loop \setbox\@picbox=\hbox{\char\count@}\ifdim\@save<\wd\@picbox
  \advance\count@-4 \repeat
 \setbox\strutbox=\hbox{\vrule height\ht\@picbox depth\dp\@picbox width\z@
   \kern\wd\@picbox}%
 \@save=.5\wd\@picbox \advance\@save-.5\linethickness
 \setbox0=\hbox to\@xdim{\ifovaltl\char\count@\else\strut\fi
  \kern-\@save\leaders\hrule height\ifovaltl\linethickness\else\z@\fi\hfil
  \leaders\hrule height\ifovaltr\linethickness\else\z@\fi\hfil\kern\@save
  \ifovaltr\advance\count@-3\char\count@\else\strut\fi\kern-\wd\@picbox}%
  \advance\count@\m@ne
 \setbox2=\hbox to\@xdim{\ifovalbl\char\count@\else\strut\fi
  \kern-\@save\leaders\hrule height\ifovalbl\linethickness\else\z@\fi\hfil
  \leaders\hrule height\ifovalbr\linethickness\else\z@\fi\hfil\kern\@save
  \ifovalbr\advance\count@\m@ne\char\count@\else\strut\fi\kern-\wd\@picbox}%
 \@save=\@ydim \advance\@save-\wd\@picbox \divide\@save 2
 \setbox\@picbox=\vbox{\box0\nointerlineskip
  \hbox to\@xdim{\vrule height\@save width\ifovaltl\linethickness\else\z@\fi
    \hfil\ifovaltr\vrule width\linethickness\kern-\linethickness\fi}%
  \nointerlineskip
  \hbox to\@xdim{\vrule height\@save width\ifovalbl\linethickness\else\z@\fi
    \hfil\ifovalbr\vrule width\linethickness\kern-\linethickness\fi}%
  \nointerlineskip\box2}%
  \@save=.5\@ydim \advance\@save.5\linethickness \leavevmode
  \kern-.5\@xdim \kern-.5\linethickness \lower\@save\box\@picbox}}

\def\cpic#1\endcpic{\vcenter{\hbox{\beginpic#1\endpic}}}


\newdimen\@xi \newdimen\@xii \newdimen\@xiii \newdimen\@xiv
\newdimen\@xpt \newdimen\@xoldpt
\newdimen\@yi \newdimen\@yii \newdimen\@yiii \newdimen\@yiv
\newdimen\@ypt \newdimen\@yoldpt
\def\squine(#1,#2,#3,#4,#5,#6){\setbox\@picbox\hbox{\tencirc q}%
 \global\@xoldpt=#1\unitlength \global\@yoldpt=#4\unitlength \kern\@xoldpt
 \@xi=\@xoldpt \@xii=#2\unitlength \@xiii=#3\unitlength
 \@yi=\@yoldpt \@yii=#5\unitlength \@yiii=#6\unitlength
 \squinerec
 \@xpt=#3\unitlength \@ypt=#6\unitlength \@addpoint
 \raise\@ypt\copy\@picbox}
\newif\iffar
\def\squinerec{\farfalse \testnear\@xi\@xiii \testnear\@yi\@yiii
 \iffar \decast \fi}
\def\testnear#1#2{\@save=#1\advance\@save-#2%
 \ifdim\@save<\z@ \@save=-\@save\fi \ifdim\@save>\p@ \fartrue \fi}
\def\decast{\@xpt=\@xi \advance\@xpt\@xii \divide\@xpt2
 \advance\@xii\@xiii \divide\@xii2
 \@xiv=\@xpt \advance\@xiv\@xii \divide\@xiv2
 \@ypt=\@yi \advance\@ypt\@yii \divide\@ypt2
 \advance\@yii\@yiii \divide\@yii2
 \@yiv=\@ypt \advance\@yiv\@yii \divide\@yiv2
 \begingroup\@xii=\@xpt \@xiii=\@xiv
  \@yii=\@ypt \@yiii=\@yiv \squinerec\endgroup
 \@xpt=\@xiv \@ypt=\@yiv \@addpoint
 \@xi=\@xiv \@yi=\@yiv \squinerec}
\def\@addpoint{
 \global\advance\@xoldpt-\@xpt \wd\@picbox=-\@xoldpt
 \raise\@yoldpt\copy\@picbox \global\@xoldpt=\@xpt \global\@yoldpt=\@ypt}

\catcode`\@=\oldatcat

\noindent
\begin{abstract}
In this paper we shall show that there exists
a polynomial unimodal map $f\colon [0,1]\to [0,1]$ which is
\begin{itemize}
\item non-renormalizable
(therefore for each $x$ from a residual set,
$\omega(x)$ is equal to an interval),
\item for which
$\omega(c)$ is a Cantor set and
\item  for which $\omega(x)=\omega(c)$ for Lebesgue almost all $x$.
\end{itemize}
So the topological and the metric attractor of such a map
do not coincide. This gives the answer to a question posed
by Milnor \cite{Mil}.
\end{abstract}

\sect{Introduction}
One of the central themes in the theory of dynamical systems
is the concept of attractors. However, there is no complete consensus
about the `correct' definition of this notion.
In particular it is not
clear whether an attractor should attract a topologically big set or
a set which is large in a metric sense.
So, if $f\colon M\to M$ is a dynamical system defined on
a manifold $M$, then
we could define a closed forward invariant set $X$ to be
a {\it topological}  respectively a {\it metric attractor} if
\begin{enumerate}
\item its basin
$$B(X)=\{x\st \omega(x)\subset X\}$$
contains a residual subset of an open neighbourhood of $X$,
respectively $B(X)$ has positive Lebesgue measure;
\item there exists no closed forward invariant
set $X'$ which is strictly included
in $X$ for which $B(X)$ and $B(X')$ coincide
up to a meager set respectively up to a set of measure zero.
\end{enumerate}
Here $\omega(x)$ is the set of limit points of $f^n(x)$ as $n\to \infty$.
Moreover, we say that $A$ is a residual (resp. meager) set if it
is the countable intersection (union)
of open dense (closed nowhere dense) sets.
For a discussion on these definitions, see \cite{Mil}.
If $X$ is a periodic attractor, a hyperbolic attractor,
a `Feigenbaum attractor' (see for example \cite{MS} and for
the invertible case see \cite{GST}), or one of the
known strange attractors, see \cite{BC}, then $X$ is both
a metric and a topological attractor. Of course, there
are some pathological cases: for example the horseshoe of
a $C^1$ diffeomorphism can have positive Lebesgue measure
and certainly is no topological attractor, see \cite{Bow}.
In this paper we present a non-pathological example
for which the distinction does matter.
More precisely, we want to show that
there exists a smooth discrete dynamical system $f\colon M\to M$
where $M$ is a smooth manifold
with an `absorbing Cantor set' $X$  (this terminology
comes from \cite{GJ}). This means that $X$ is a closed
forward invariant minimal set $X\subset M$
with zero Lebesgue measure, such that its basin $B(X)$
has positive Lebesgue measure but its complement is
a residual set. As far as we know this example is the first
smooth dynamical system with such an `absorbing Cantor set'.

In our case $M=[0,1]$ and $f$ is a smooth unimodal interval map --
this means $f$ has one extremal point --
and for simplicity we shall also assume that  $f(0)=f(1)=0$.
A prototype of such map is
$$f(x)=\lambda\left[ 1 - |2x-1|^\ell\right]$$
where $\lambda>0$ is chosen so that
$f$ maps the interval $[0,1]$ inside itself and
$f$ has the so-called Fibonacci-type dynamics. We shall define
this in the next section.

There are many publications in which it was conjectured
that a smooth map $f\colon [0,1]\to [0,1]$
cannot have an absorbing Cantor set. (We should note, however, that in 1992 Misha Lyubich and
Folkert Tangerman made computer estimates suggesting that 
absorbing Cantor sets do exist for Fibonacci maps of the form $x\mapsto x^6+c_1$.)
Moreover,  there are several results which prove that
these sets cannot exist in particular cases,
see \cite{JS1}, \cite{LM} and in the general quadratic unimodal case
\cite{L1} when $\ell=2$. We shall show that
absorbing Cantor sets do exist when $\ell$ is a large real number.
\bigskip

\noindent
{\bf Main Theorem}\\
{\em There exists $\ell_0$ with the following property.
Let $f$ be a $C^2$ unimodal interval map
with a critical point of order $\ell\ge \ell_0$
and with the Fibonacci combinatorics.
Then $f$ has an absorbing Cantor attractor $X$. }
\bigskip

Here we say that $c$ is a critical point of a $C^2$ map $f$
if $Df(c)=0$ and the {\it order of the critical point} is said to be
$\ell$ if there exists a $C^2$ diffeomorphism $\phi$ between
two neighbourhoods of $c$ such that
$$f\circ \phi(x)=f(c)-|x-c|^\ell$$
for $x$ close to $c$. It is easy to show that
our methods also give examples of multimodal smooth interval maps
for which each critical point is quadratic and
which have an absorbing Cantor set: simply choose the
map so that the return map near some critical point is
a unimodal map of Fibonacci-type
while the orbit of this critical point
contains at least $\ell$ other critical points.
However, it is not clear whether absorbing Cantor attractors also
appear generically in one-parameter families:
\bigskip

\noindent
{\bf Question} \quad
{\em Does the space of smooth maps $f\colon [0,1]\to [0,1]$
with an absorbing Cantor set form a codimension-one subset
of the space of all smooth interval maps?}
\bigskip

Of course, it follows from the Main Theorem
that there exists on each
smooth manifold a smooth mapping with an absorbing
Cantor set. We conjecture that one can also construct
invertible examples:
\bigskip

\noindent
{\bf Conjecture} \quad
{\em For each $n\ge 2$ dimensional smooth manifold $M$,
there exists a diffeomorphism $f\colon M\to M$
which has an absorbing Cantor set.}
\bigskip

In the complex one-dimensional direction there
are related results:
\bigskip

\noindent
{\bf Theorem}\quad \cite{NS2}\\
{\em For each sufficiently large even integer $\ell$
there exists $c_1\in \rz$ such that the map
$f(z)=z^\ell+c_1$ has the following properties:
\begin{itemize}
\item the set $\omega(0)$ is a Cantor set with zero Lebesgue measure;
\item the set of points $z\in \cz$ for which
$\omega(z)$ is contained in $\omega(0)$ has positive Lebesgue measure;
\item the set of points whose forward iterates remain bounded
has no interior.
\end{itemize}
In particular, the Julia set of
$z\mapsto z^\ell+c_1$ has positive Lebesgue measure.
This map has the Fibonacci dynamics (to be defined in the next section).}

\bigskip

\subsection{Some comments on the Main Theorem and its proof}
In fact, the attractor $X$ from the Main Theorem is equal to
$\omega(c)$ and this set
has zero Lebesgue measure, see \cite{Mar} and also \cite{MS}.
If the map $f$ from the
Main Theorem is a unimodal polynomial with a unique
critical point in $\cz$ (or if has negative Schwarzian
derivative and $f$ has no attracting fixed points)
then $B(X)$ has full Lebesgue measure and
its complement is a residual set.
We should remark that a smooth map as above may have
one or more periodic attractors, but that even then
the attractor $X$ has a basin
which attracts a set of positive Lebesgue measure
(and the critical point is density point of $B(X)$).
This is not completely surprising because
$\omega(c)$ is not accumulated by periodic attractors,
see \cite{MMS} and also \cite{MS}[Chapter IV].

In the theory of unimodal interval maps with
negative Schwarzian derivative of $f$, i.e., with
$$Sf(x)=\frac{D^3f(x)}{Df(x)}-\frac{3}{2}\frac{D^2f(x)}{Df(x)}<0$$
and for which the order of the critical point is finite,
one has a well-known classification,
see \cite{G}, \cite{BL}, \cite{keller1} and also \cite{MS}.
\begin{enumerate}
\item $f$ has a stable periodic orbit $O$
which is both a topological and metric attractor;
\item $f$ is infinitely renormalizable, i.e., there
exists a nested sequence of intervals $I_n\ni c$
shrinking to $c$ and a sequence of integers $q(n)\to \infty$
such that $I_n,\dots,f^{q(n)-1}(I_n)$ are disjoint
and $f^{q(n)}(I_n)\subset I_n$. In this case
$\omega(c)$ is a Cantor set of zero Lebesgue measure
which is both a topological and metric attractor;
\item $f$ is not infinitely renormalizable.
In this case there exists a cycle of intervals $Z$
(a finite union of intervals) such that
$B(Z)$ is dense and has full Lebesgue measure.
The set $Z$ is a topological attractor,
but not necessarily a metric attractor:
in principle, there could be a Cantor set $X\subset Z$ such that
$B(X)$ has full Lebesgue measure (but is not dense).
\end{enumerate}
From our theorem it follows that the possibility mentioned in
the last case really does occur if $\ell$ is large.
In the quadratic case, i.e.
$\ell=2$, the results of \cite{L1} imply that
$Z$ is a metric attractor as well.

Any map with an absorbing Cantor set has no absolutely
continuous invariant probability measure, because
Lebesgue almost all points wander densely on the support of the
measure by the Birkhoff Ergodic Theorem.
If the Schwarzian derivative of $f$
is negative and $\ell=2$ then it is shown in \cite{LM}
that $f$ has an absolutely continuous invariant probability measure
by showing that the summability condition from \cite{NS1}
is satisfied. In particular, $f$ has no absorbing Cantor set
in this case. The methods of proof in \cite{LM} are a mixture
of real tools and tools from the theory of complex analysis
and hyperbolic geometry.
This result was generalized in \cite{fibo}:
in that paper it was shown that the same results hold for
$1<\ell\le 2+\epsilon$
provided $\epsilon>0$ is small. The tools in \cite{fibo} are entirely
based on real estimates, and also no use is made of \cite{NS1}
(because the summability condition fails if $\ell>2$).

As mentioned, our result implies that $f$ has no
absolutely continuous invariant probability measure
for $\ell$ large. In fact, as Henk Bruin has shown in \cite{Br},
this already follows from Proposition~\ref{onestepbounds}.

We expect that the methods of this paper can be extended
to show that for Fibonacci maps of `bounded type'
(a notion which we shall discuss in the section about
the combinatorial properties of Fibonacci maps)
with a rather flat critical point, the same result holds.

Let us now give an outline of the proof that $\omega(x)$
is equal to the Cantor set $\omega(c)$ for Lebesgue almost all $x$.
\begin{itemize}
\item
First we will show that there exists
a nested sequence of intervals $(u_n,\hat u_n)$
containing $c$ and that the size of the annulus
$A_n=(u_n,\hat u_n)\setminus (u_{n+1},\hat u_{n+1})$
is very small compared to the size of $(u_{n+1},\hat u_{n+1})$
if the order $\ell$ of the critical point is large.
\item Next we let $I_n,\hat I_n$ be the components of $A_n$ and
show that some iterate $f^{S_n}$ of $f$ maps $I_n$
diffeomorphically inside $\cup_{k\ge n-2}(I_k\cup \hat I_k)$
and that this map is not `too' non-linear.
Because of 1) this implies that `most' points are mapped
closer to $c$ by this iterate.
\item Finally, we combine 1), 2) and a kind of random walk argument
to show that typical points are in the basin of $\omega(c)$.
\end{itemize}

\sect{Combinatorial properties of the Fibonacci map}
\label{sectop}
In this section we shall define and state some properties of
the Fibonacci map. It is well-known that maps with these
properties exist, see \cite{HK} or
\cite{LM} and also the sequel to this paper. In the companion paper
\cite{NS2} we shall construct such a map `by hand'.
Let $f\colon [0,1]\to [0,1]$ be a unimodal map with $f(0)=f(1)=0$.
For each $x\ne c$ there exists a `symmetric' point $\hat x\ne x$
with $f(\hat x)=f(x)$. For $i\ge 0$ and $x\in [0,1]$,
let $x_i=f^i(x)$ and
choose $x_{-i}\in f^{-i}(x)$ so that
the interval connecting this point to $c$ contains no other points
in the set $f^{-i}(x)$. Note that if $c$ is not a periodic point
there are always precisely two such points $c_{-i}$
(which are symmetric with respect to each other).
Let $S_0=1$ and define $S_i$ inductively by
$$S_i=\min\{k\ge S_{i-1};\, c_{-k}\in (c_{-S_{i-1}},\hat c_{-S_{i-1}})\}.$$
$f$ is called a {\it Fibonacci} map
if the sequence $S_i$ coincides with
the Fibonacci numbers: $S_0=1$, $S_1=2$ and $S_{k+1}=S_k+S_{k-1}$, i.e.,
the sequence $1,2,3,5,8,\dots$.
The proof of the following proposition can be found in
\cite{LM}, \cite{fibo} and also in \cite{NS2}.

Let us denote by $z_k$ the nearest point to $c$
in the set $f^{-S_k}(c)$. It should be clear from the context whether
$z_k$ is to the left or right of $c$.
Moreover, for $x\in [0,1]$ let us write
$$x^f=f(x)$$
(usually, $x$ will be close to $c$ and so $x^f$ close
to $c^f=f(c)$).

\begin{prop}\label{top1}
A Fibonacci map $f\colon [0,1]\to [0,1]$
satisfies the following properties.

\begin{itemize}
\item $f$ is non-renormalizable;
\item $c_{S_k}$ and $c_{S_{k+2}}$ are on opposite sides
of $c$.
\item $c_{S_n}\in (c_{S_{n-1}},\hat c_{S_{n-1}})$
and $c_i \notin (c_{S_{n-1}},\hat c_{S_{n-1}})$
for each $0<i<S_n$.
\item
$c_{-S_n}\in (c_{-S_{n-1}},\hat c_{-S_{n-1}})$
and $c_{-i} \notin (c_{-S_{n-1}},\hat c_{-S_{n-1}})$
for each $0<i<S_n$.
\item If $T$ is the maximal interval adjacent
to $c$ such that $f^{S_k}_{|T}$ is monotone,
then $f^{S_k}(T)=(c_{S_k},c_{S_{k-2}})$.
\item If $T_k \owns c^f$ is the largest interval
on which $f^{S_k-1}_{|T_k}$ is monotone, then
$$T_k=(z_{k-1}^f,t_k^f)$$
where $t_k^f>c^f$ and
$f^{S_k-1}(T_k)=(c_{S_{k-2}},c_{S_{k-4}})$
(note that $t_k^f$ is not the $f$-image of some point
$t_k$, so this notation is just to suggest that
$t_k^f$ is close to $c^f$).
\item $T_k,\dots,f^{S_k-1}(T_k)$ has intersection multiplicity
$\intersect$ (this means that each point of $[0,1]$ is contained in
at most $\intersect$ of these intervals).
\end{itemize}
\end{prop}
\pr
For the proof of this result we refer the reader to
\cite{fibo} and \cite{LM}. The statement about the
disjointness can be found in \cite{LM}[Lemma 4.3].
\qed

From the fact that $c_{-1}$ exists it follows that
$f$ has a orientation reversing fixed point $q$.
Let us define inductively a sequence of points $u_n$ as follows.
Let $u_0=q$ and let us define $u_{n+1}$ to be nearest point
to $c$ with
$$u_{n+1}\in f^{-S_n}(u_n)$$
so that $u_{n+1}$ is on the same side of $c$ as $c_{S_{n+1}}$.
In particular, $u_1=\hat u_0=\hat q$.
Moreover, let $\tilde u_{k+1}$ be the point
in $\{u_{k+1},\hat u_{k+1}\}$ which is on the same side
of $c$ as $u_k$.

Furthermore, let
$$
y_n=\fSn{n}(\cSn{n+2})\ ,\quad y_n^f=f(y_n).
$$

\begin{prop}\label{top2} A Fibonacci map $f\colon [0,1]\to [0,1]$
satisfies the following properties.

\begin{itemize}
\item $f^{S_n}(u_{n+1})=u_n$ and $f^{S_n}(u_n)=u_{n-2}$;
\item in particular, $f^{S_n}$ maps $(\tilde u_{n+1},u_n)$
diffeomorphically onto $(u_n,u_{n-2})$ (note that this last interval
contains $c$);
\item the points $u_n^f$, $c_{S_n}^f$, $c_{S_n+S_{n+2}}^f$, $y_n^f$ and $z_n^f$
are ordered as in the picture below
(we state the ordering near $c^f$ rather than near $c$
so that we do not need to be careful about on which side of $c$
these points lie).
\end{itemize}
\end{prop}
\hbox to \hsize{\hss\unitlength=5mm
\beginpic(20,4)(0,0) \let\ts\textstyle
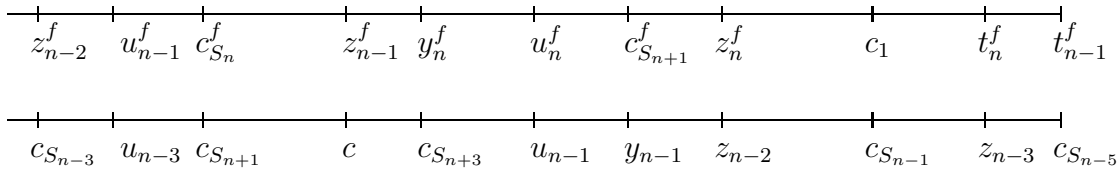
\put(3,-3.3){{\it Figure 2.1: Points and their images under $f^{S_{n-1}-1}$}.}
\put(-4,3){\line(1,0){28}}
\put(-3.2,2.8){\line(0,1){0.4}} \put(-3.4,2){$z_{n-2}^f$}
\put(-1.2,2.8){\line(0,1){0.4}} \put(-1,2){$u_{n-1}^f$}
\put(1.2,2.8){\line(0,1){0.4}} \put(1,2){$c_{S_{n}}^f$}
\put(5,2.8){\line(0,1){0.4}} \put(4.9,2){$z_{n-1}^f$}
\put(7,2.8){\line(0,1){0.4}} \put(6.9,2){$y_n^f$}
\put(10,2.8){\line(0,1){0.4}} \put(9.9,2){$u_n^f$}
\put(12.5,2.8){\line(0,1){0.4}} \put(12.4,2){$c_{S_{n+1}}^f$}
\put(15,2.8){\line(0,1){0.4}} \put(14.8,2){$z_{n}^f$}
\put(19,2.8){\line(0,1){0.4}} \put(18.8,2){$c_1$}
\put(22,2.8){\line(0,1){0.4}} \put(21.8,2){$t_n^f$}
\put(24,2.8){\line(0,1){0.4}} \put(23.8,2){$t_{n-1}^f$}

\put(-4,0.2){\line(1,0){28}}
\put(-3.2,0){\line(0,1){0.4}} \put(-3.4,-0.8){$c_{S_{n-3}}$}
\put(-1.2,0){\line(0,1){0.4}} \put(-1,-0.8){$u_{n-3}$}
\put(1.2,0){\line(0,1){0.4}} \put(1,-0.8){$c_{S_{n+1}}$}
\put(5,0){\line(0,1){0.4}} \put(4.9,-0.8){$c$}
\put(7,0){\line(0,1){0.4}} \put(6.9,-0.8){$c_{S_{n+3}}$}
\put(10,0){\line(0,1){0.4}} \put(9.9,-0.8){$u_{n-1}$}
\put(12.5,0){\line(0,1){0.4}} \put(12.4,-0.8){$y_{n-1}$}
\put(15,0){\line(0,1){0.4}} \put(14.8,-0.8){$z_{n-2}$}
\put(19,0){\line(0,1){0.4}} \put(18.8,-0.8){$c_{S_{n-1}}$}
\put(22,0){\line(0,1){0.4}} \put(21.8,-0.8){$z_{n-3}$}
\put(24,0){\line(0,1){0.4}} \put(23.8,-0.8){$c_{S_{n-5}}$}
\endpic\hss}

\vskip 2.2cm
\pr
The proof of these statements can be found in \cite{fibo}.
It can be derived from Figure 2.2 below.
\qed

\vskip -1cm

\hbox to \hsize{\hss\unitlength=1.5mm
\beginpic(70,110)(0,0) \let\ts\textstyle
\put(-5,40){\line(1,0){60}}
\put(30,1){\line(0,1){93}}
\put(25,-4){{\it Figure 2.2.}}
\put(25,48){\squine(0, -20, -30, 0, -50, -43)}
\put(25,48){\squine(0, 20, 30, 0, 50, 43)}
\put(32,75){$f^{S_{n-1}-1}$}

  \put(-2,40){\line(0,-1){36}} \put(-1,37.5){$z_{n-2}^f$}
\put(8,40){\line(0,-1){26.8}} \put(7.5,42){$u_{n-1}^f$}
  \put(17,40){\line(0,-1){11}} \put(13.7,37.5){$c_{S_n}^f$}
  \put(16.3,42){$z_{n-1}^f$}
  \put(25,40){\line(0,1){8}} \put(24,37.5){$y_n^f$}
  \put(27.5,40){\line(0,1){14}} \put(26.5,37.5){$c_{S_{n+1}}^f$}
  \put(31,38){$c^f$}
  \put(52,40){\line(0,1){52}} \put(51,37.5){$t^f$}

  \put(30,4){\line(-1,0){32}} \put(31,4){$c_{S_{n-3}}$}
\put(30,13.2){\line(-1,0){22}} \put(31,13.2){$u_{n-3}$}

  \put(30,29){\line(-1,0){13}} \put(31,29){$c_{S_{n+1}}$}
  \put(30,48){\line(-1,0){5}} \put(31,48){$c_{S_{n+3}}$}
  \put(30,54){\line(-1,0){2.5}} \put(31,54){$y_{n-1}$}
  \put(31,60){$c_{S_{n-1}}$}
  \put(30,92){\line(1,0){22}} \put(31,90){$c_{S_{n-5}}$}
\endpic\hss}
\bigskip\bigskip

If $f$ has no wandering intervals
(and this is the case under the present assumptions, see
\cite{MS}[Chapter IV]), then $\omega(c)$ is a minimal Cantor set.

\sect{The estimates}
In this section we shall estimate the rate of approach of the
sequences $u_k^f, c_{S_k}^f, z_k^f$ to $c^f$.
The basic tool is that of the distortion of cross-ratios.

\begin{remark}
Since $f$ is non-flat at $x=c$ we
can assume (by applying  a suitable $C^2$ coordinate
change) that $f$ is of the form
$$f(x)=f(c)-|x-c|^\ell$$
near $x=c$. This will simplify some of the estimates somewhat.
\end{remark}

Hence
$$\frac{|f(x)-f(c)|}{|x-c|}=M(x)$$
where $M(x)$ is a continuous function
which is equal to $|x-c|^{\ell-1}$ near $x=c$.
Moreover,
$$\frac{Df(x)}{\ell |f(x)-f(c)|/|x-c|}=1$$
near $x=c$.
We shall use these facts repeatedly.

\subsection{The cross-ratio and the Koebe Principle}
Let $j\subset t$ be intervals and let $l,r$ be the components
of $t\setminus j$. Then
the cross-ratio of this pair of intervals is defined
as
$$C(t,j):=\frac{|t|}{|l|}\frac{|j|}{|r|}.$$
Let $f$ be a smooth function  mapping
$t,l,j,r$ onto $T,L,J,R$ diffeomorphically.
Define 
\[
B(f,t,j)=\frac{|T|\, |J|}{|t|\, |j|}\,\frac{|l|\, |r|}{|L|\, |R|}\
=\frac{C(T,J)}{C(t,j)}
.
\]
It is well known that if $Sf=f'''/f'-3(f''/f')^2/2\le 0$ then $B(f,t,j)\ge 1$.
In the next proposition it is stated that this ratio
also cannot be decreased too much by a $C^2$ map $f$
with non-flat critical points.

\begin{prop}\label{crrat}
Let $f$ be a $C^2$ map with non-flat critical points.
Then there exists a function $o(\epsilon)>0$ with
$o(\epsilon)\to 0$ as $\epsilon\to 0$
such that for any intervals $j\subset t$ and any $n$
for which $f^n|t$ is a diffeomorphism one has
the following. Let $l,r$ be as above and let $L,J,R,T$
be the images of $l,j,r,t$ under $f^n$.
Then
$$B(f^n,t,j)=\frac{|T|}{|L|}\frac{|J|}{|R|}
\frac{|l|}{|t|}\frac{|r|}{|j|}\ge
\exp\left(-o(\epsilon)\cdot \sum_{i=0}^{n-1}|f^i(t)| \right)
$$
where $\epsilon=\max_{i=0}^n|f^i(t)|$.
(If $Sf<0$ then $B(f^n,t,j)>1$.)
\end{prop}
\pr See Theorem IV.2.1 in \cite{MS}.
\qed

From this it follows in particular that if  $f^n|T$ is a diffeomorphism
and $j$ is reduced to the point $x$ then
$$\frac{|Df(x)|}{|L|/|l|}\ge
\exp\left(-o(\epsilon)\cdot \sum_{i=0}^{n-1}|f^i(t)| \right)
|R|/|T|$$
and if $l$ is reduced to the point  $y$ then
$$\frac{|Df(y)|}{|J|/|j|}\le
\exp\left(o(\epsilon)\cdot \sum_{i=0}^{n-1}|f^i(t)| \right)
|T|/|R|.$$

We shall also need the following lemma.
In fact, instead of this lemma one could use the
Koebe Principle stated below (and the Koebe Principle
can be derived from the next lemma).

\begin{lemma}\label{doubc}
Let $f$ be a $C^2$ map with non-flat critical points.
Then there exists a function $o(\epsilon)>0$ with
$o(\epsilon)\to 0$ as $\epsilon\to 0$
such that for any $n$ and any interval $t$ for which
$f^n|t$ is a diffeomorphism one has the following property.
Let $j_1,j_2\subset t$ be two disjoint intervals
and let $l_i,r_i$ be the components of $t\setminus j_i$
where we assume $l_1\subset l_2$ and
$r_1\supset r_2$.
(So this is the case if $j_1$ lies to the left
of $j_2$ and $l_i$ is to the left of $j_i$ for $i=1,2$.)
Let $L_i,J_i,R_i,T$ be the images of $l_i,j_i,r_i,t$ under $f^n$.
Then
$$\frac{|j_1|}{|J_1|}
\frac{|J_2|}{|j_2|}
\ge
\O \cdot
\frac{|J_2\cup R_2||R_2|}{|J_1\cup R_1||R_1|}=
\frac{|J_2\cup R_2||R_2|}{|J_1\cup J\cup J_2\cup R_2||J\cup J_2\cup R_2|}
$$
and
$$\frac{|j_1|}{|J_1|}
\frac{|J_2|}{|j_2|}
\le
\O \cdot
\frac{|L_2||L_2\cup J_2|}{|L_1||L_1\cup J_1|}.
$$
where
$$\O =\exp\left(\pm o(\epsilon)\cdot \sum_{i=0}^{n-1}|f^i(t)| \right)
$$
and $\epsilon=\max_{i=0,\dots,n-1}|f^i(t)|$.
(If $Sf<0$ then we can take $\O=1$.)
\end{lemma}
\pr
Let $j$ be the interval connecting $j_1$ and $j_2$.
Multiplying the following two cross-ratio
inequalities from the previous proposition,
the result follows immediately.
$$
\frac{|J|}{|J_1|}
\frac{|J_1\cup J\cup J_2\cup R_2|}{|J_2\cup R_2|}
\ge
\O
\frac{|j|}{|j_1|}
\frac{|j_1\cup j\cup j_2\cup r_2|}{|j_2\cup r_2|}
$$
and
$$
\frac{|J_2|}{|J|}
\frac{|J\cup J_2\cup R_2|}{|R_2|}
\ge \O
\frac{|j_2|}{|j|}
\frac{|j\cup j_2\cup r_2|}{|r_2|}.
$$
The second inequality follows similarly.
\qed


\vskip 0.5cm
\hbox to \hsize{\hss\unitlength=4.5mm
\beginpic(19,10)(0,0) \let\ts\textstyle
\put(3,1){{\it Figure 3.1: Intervals $j_i,l_i,r_i$ and their images.}}
\put(0,4){\line(1,0){19}}
\put(2,4){\line(0,1){0.2}} \put(3.3,4.5){$L_1$}
\put(5,4){\line(0,1){0.2}} \put(6.3,4.5){$J_1$}
\put(8,4){\line(0,1){0.2}} \put(9.3,4.5){$J$}
\put(11,4){\line(0,1){0.2}} \put(12.3,4.5){$J_2$}
\put(14,4){\line(0,1){0.2}} \put(15.3,4.5){$R_2$}
\put(17,4){\line(0,1){0.2}}
\put(0,8){\line(1,0){19}}
\put(2,8){\line(0,1){0.2}} \put(3.3,8.5){$l_1$}
\put(5,8){\line(0,1){0.2}} \put(6.3,8.5){$j_1$}
\put(8,8){\line(0,1){0.2}} \put(9.3,8.5){$j$}
\put(11,8){\line(0,1){0.2}} \put(12.3,8.5){$j_2$}
\put(14,8){\line(0,1){0.2}} \put(15.3,8.5){$r_2$}
\put(17,8){\line(0,1){0.2}}
\endpic\hss}
\vskip 0.5cm

We should remark that if we take $T_k$ to be the maximal interval
containing $c_1$ on which $f^{S_k-1}$ is a diffeomorphism,
then from  Proposition~\ref{top1},
$$\sum_{i=0}^{S_k-1}|f^i(T_k)|\le \intersect.$$
So this implies that we can apply the previous results immediately
to $f^{S_k-1}|T_k$.
In fact, $\max_{i=0}^{S_k-1}|f^i(T_k)|\to 0$
as $k\to \infty$:
\bigskip
\begin{lemma}\label{contrpr}
For each $\epsilon>0$ there exists
$\delta>0$ such that if
$f^n(I)$ is not contained in an immediate basin of
a periodic attractor and $|f^n(I)|\le \delta$,
then $\max_{i=0}^{n-1}|f^i(I)|\le \epsilon$.
\end{lemma}
\pr If the lemma is not satisfied, then
there exists a sequence of intervals $I_i$ with
$|I_i|\ge \epsilon$ and a sequence $n(i)$ with
$|f^{n(i)}(I_i)|\to 0$ where $f^{n(i)}(I_i)$ is not completely contained
in the immediate basin of some periodic attractor. By taking
subsequences, there exists an interval $I$ such that
$\inf_{i\ge 0} |f^i(I)|=0$ and such that $I$ is not completely contained
in the basin of a periodic attractor. This is impossible
because $f$ has no wandering intervals,
see \cite{MS}[Chapter IV, Theorem A].
Indeed, by the Contraction Principle, see  \cite{MS}[IV.5.1]
if $I$ is an interval with $\inf_{i\ge 0}|f^i(I)|=0$ then either
$I$ is completely contained in the basin of a periodic attractor or
a wandering interval.
\qed

In fact, we shall also have to estimate the cross-ratio
distortion of iterates of $f$ which are not of the form $f^{S_i}$.
For this we shall need the Koebe Principle and an estimate
on the total size of orbits of some intervals.
Let us say that an interval $T$ contains 
a $\tau$ -{\it scaled neighbourhood} of an interval $J\subset T$
if each component of $T\setminus J$ has at least size
$\tau|J|$.

\bigskip
\begin{prop} [Koebe Principle]
Let $f$ be a $C^2$ map with non-flat critical points.
Then there exists a function $o(\epsilon)>0$ with
$o(\epsilon)\to 0$ as $\epsilon\to 0$
such that for any intervals $j\subset t$ and any $n$
for which $f^n|t$ is a diffeomorphism one has
the following.
If $f^n(t)$ contains a $\tau$-scaled
neighbourhood of $f^n(j)$
then
\beq
\label{koebee}
\frac{|Df^n(x)|}{|Df^n(y)|}\le
\left[\frac{1+\tau}{\tau}\right]^2
\exp\left( o(\epsilon)\cdot \sum_{i=0}^{n-1}|f^i(j)| \right)
\eeq
for each $x,y\in J$
where $\epsilon=\max_{i=0}^n|f^i(t)|$.

Moreover, if $f^n(t)$ contains no periodic attractor
then there exists $t'$ with $j\subset t'\subset t$
for which
$$f^n(t')\text{ is a }\tau/2-\text{scaled neigbourhood of }f^n(j)$$
and
$$|f^i(t')|\le K|f^i(j)|
$$
for all $i=0,1,\dots,n$.
Here $K$ depends on $f$, $\tau$, $\epsilon$ and $\sum_{i=0}^{n-1}|f^i(j)|$.
\end{prop}
\pr This is a combined statement of
Theorem IV.3.1 and Theorem IV.1.1 in \cite{MS} and of
Lemma 8.3 from \cite{Str}. (Note that we do not assume 
that $\sum |f^i(t)|$ is bounded but merely that 
$\sum |f^i(j)|$ is bounded.
\qed

In the next proposition we shall give a condition
for orbits of intervals to have a finite total length.
We shall need this proposition only in the case
that the Schwarzian derivative of $f$ is not negative
to estimate the term $\O$ in Lemma \ref{doubc}.

\bigskip

\begin{prop}\label{sumlength}
Let $f$ be a $C^2$ map with non-flat critical points.
Then for each $\tau,S>0$
there exist constants $\delta,S'>0$
such that the following holds.
Let $x$ be a recurrent point of $f$, let
$U$ be an interval neighbourhood around $x$ of size $<\delta$
and $D\subset U$ be some disjoint union of intervals $I_i$.
Let $F$ be a map defined on $D=\cup I_i$ such that
for each interval $I_i$
there exists an integer $j(i)$ and an interval $T_i\supset I_i$
such that
\begin{enumerate}
\item 
$F|I_i=f^{j(i)}$ maps $I_i$ onto some union of intervals $I_i$
and $F(I_i)$ contains at least two of those intervals;
\item 
$f^{j(i)}|T_i$ is a diffeomorphism
and $I_k\subset f^{j(i)}(I_i)$ implies $T_k\subset f^{j(i)}(T_i)$;
\item 
$f^{j(i)}(T_i)$ contains a $\tau$-scaled neighbourhood of
each interval $I_k\subset f^{j(i)}(I_i)$;
\item 
for each interval $I_j\subset F(I_i)$ one has
$|I_j|\le (1-\frac{1}{S})\cdot |F(I_i)|$;
\item 
$\sum_{m=0}^{j(i)-1}|f^m(T_i)|\le S$.
\end{enumerate}
Then for each $n\in \nz$ and each component $J$ of the domain
of $F^n$ one has $F^n|J=f^j$ for some $j\in \nz$
and there exists an interval $T'\supset J$
for which $F^n(T')$ contains a $\tau/2$-scaled neighbourhood of
each element $I_k\subset F^n(J)$
and
$$\sum_{m=0}^{j-1}|f^m(T')|\le S'.$$
\end{prop}
\pr
The idea of the proof of this proposition is essentially the same
as that of \cite{Str}. The proof of this proposition
can be substantially simplified if $f$ has
negative Schwarzian derivative: in this case it is not necessary to
choose $\delta$ small. (In fact we do not even
need this lemma in that case.)
However, in the general case, $f$ could
for example have a periodic interval (corresponding to basins
of periodic attractors). This complicates matters to some
extend.

Fix $\tau$ and $K$. Since $f$ is $C^2$
each periodic point $p$ of $f$ of sufficiently large
period $k$ is repelling,
see \cite{MS}[Theorem IV.B].
In particular, this holds for all periodic points
which are in a $\delta$ neighbourhood of $x$, provided
$\delta>0$ is sufficiently small. For this reason we
shall be able to apply Lemma~\ref{contrpr}.

Now let $\I$ be the partition of the domain of $F$
of the intervals $I_i$ and define inductively $\I_0=\I$ and
$$\I_n=\I_0\vee F^{-1}\I_0 \vee \dots \vee F^{-n}\I_0.$$
So each element $J$ of $\I_n$ is an interval
which $F^n$  maps diffeomorphically onto some interval $I_i$
and each of these intervals is contained in $U$
where $|U|<\delta$.
Because of properties 2) and 3) there exists
$j\in \nz$ with $F^n|J=f^j$
and an interval $T\supset J$ which is mapped diffeomorphically
onto a $\tau$-scaled neighbourhood of $f^j(J)=I_k\in \I_0$.
It suffices to show that there exists
$S'$ such that
\beq
\label{suff}
\sum_{i=0}^j|f^i(J)|\le S'\text{ for each }J\in \I_n.
\eeq
Indeed, the components of the domain of $F^n$ are elements from
$\I_{n-1}$ and are mapped by $F^{n-1}$ into an of elements of $\I_0$.
However, because of property 5) the length of the remaining intervals
up to the $F^n$-th iterate have uniformly bounded sum.

First we claim that there exists $\kappa<1$ such that
\beq
\label{ref}
\text{if $J\in \I_1$ is contained in $I_i\in \I_0$
then $|J|\le \kappa|I_i|$.}
\eeq
This holds since $F(J)$ is equal to an interval $I_k\in \I_0$ while
properties 3) and 4) imply that there exists an interval $J'$ with
$J\subset J'\subset I_i$ for which
i) $F(J')$ is contained inside a $\tau/2$-scaled neighbourhood of $F(J)=I_k$
and ii) a definite proportion of $F(J')$ is outside
$F(J)=I_k$. Moreover, because of 5) and the Koebe Principle there exists
a universal constant $K_0<\infty$ such that
$$\sup_{x,y\in J'}\frac{|DF(x)|}{|DF(y)|}\le K_0.$$
Combining this proves (\ref{ref}).

By using a `telescope argument' we can improve this statement
and show by induction that
there exists $\kappa<1$ such that
for each $n\in \nz$ there exists $\delta>0$ such that
if $|U|<\delta$,
$J\in \I_n$ and $J$ is contained in $I_i\in \I_0$ then
\beq
\label{taun}
|J|\le \kappa^n|I_i|.
\eeq
For $n=0$ there is nothing to prove.
So assume the statement holds for $n-1$
and consider $J\in \I_n$.
If $F^n|J=f^j$ and $T\supset J$ so that $f^j|T$ is a diffeomorphism and
$f^j(T)$ is a $\tau$-scaled neighbourhood of $f^j(I)=I_k\in \I_0$ then
\beq
\label{ww}
\sum_{i=0}^{j-1}|f^i(J)|\le  n\cdot S
\text{ and }\max_{i=0,\dots,j-1}|f^i(T)|=
o(|f^j(T)|),
\eeq
where $o(t)$ is a function so that $o(t)\to 0$ if $t\downarrow 0$.
Here we have used respectively
property 5) and the previous Lemma \ref{contrpr}.
(We should note that $f^n(T)\subset U$ and so
$|f^j(T)|=|F^n(J)|\le |U|\le \delta$.)
Hence, by the Koebe Principle,
there exists $K_1$ (which only depends on $\tau$)
such that for each $n$ 
\beq
\label{koe2}
\frac{|DF^n(x)|}{|DF^n(y)|}\le K_1
\eeq
for all $x,y\in J$ provided $\delta$ (and hence
$F^n(T)\subset U$) is sufficiently small.
(To get $K_1$ uniform
we shrink $\delta$ for increasing $n$; by (\ref{ww}) and
(\ref{koebee}) this avoids the constants in the Koebe Principle to grow.)
Now $F^{n-1}$ maps each element of $\I_{n-1}$ diffeomorphically onto
some element of $\I_0$  and each element of $\I_{n}$ onto 
an element of $\I_1$. From this, (\ref{koe2}) and (\ref{ref})
it follows that each element $J$ of $\I_n$
is a definite factor smaller than the element $I\in \I_{n-1}$ containing
$J$. This proves (\ref{taun}).

Now of course (\ref{taun}) does not suffice because
$\delta$ (and therefore the size of $U$) depends on $n$.
Therefore,
let us fix $n_0$ so large that
$$\kappa^{-n_0}\ge 4K_2\text{ where }K_2=\left[ \frac{1+\tau}{\tau}\right]^2$$
and write $G=F^{n_0}$. If $J$ is an element of
$\I_{kn_0}$ and $G^i(J)\supset J$ for some
$0\le i \le k$ then
\beq
\label{twof}
|DG^i(x)|\ge 2\text{ for all }x\in J.
\eeq
Indeed, we may assume that $i$ is minimal
and then $J,\dots,G^{(i-1)}(J)$ are disjoint.
If $G^i=f^j$ then this gives that
$J,\dots,f^j(J)$ have intersection multiplicity
bounded by $n_0$. (This means that each point is contained
in at most $n_0$ of these intervals.)
Therefore, and since $f^j$ maps some interval $T\supset J$
onto a $\tau$-scaled neighbourhood of $f^j(J)$, it follows from
the Koebe Principle that
\beqa
|DG^i(x)|
&\ge&
\exp\left(-o(\epsilon)\cdot \sum_{i=0}^{j-1} |f^i(J)|\right)
\frac{1}{K_2}\frac{|G^i(J)|}{|J|}\\
&\ge&
\exp\left(-o(\epsilon)\cdot n_0 \right)
\frac{1}{K_2}\frac{|G^i(J)|}{|J|}\ge
\frac{1}{2}\frac{1}{K_2}\kappa^{-n_0}\ge 2
\eeqa
for each $x\in J$ provided $|f^j(J)|=|G^i(J)\le |U|\le \delta$
is sufficiently small. (This last inequality implies
that $\epsilon=\max|f^i(J)|$ is small when $\delta$
is small.) Hence, if some interval returns then its size
has increased by a uniform factor; as we shall now show this implies
the total length of the intervals remains bounded.
Indeed, consider again $J\in \I_{kn_0}$.
Then
\beq
\label{satel}
\sum_{i=0}^{k-1}|G^i(J)|\le \frac{1}{1-1/2}=2.
\eeq
This is because $G^{i_1}(J)\cap G^{i_2}(J)\ne \emptyset$
with $i_1<i_2\le k$ implies that
$G^{i_1}(J)\subset G^{i_2}(J)$.
Moreover, if $G^{i_1}(J),G^{i_2}(J)\subset G^{i_3}(J)$
and $i_1<i_2\le i_3\le k$, then there exists
$J'\supset G^{i_1}$ (which is an interval from a partition
of the form $\I_{hn_0}$ with $h\in \{0,1,\dots,k\}$)
such that $G^{i_2-i_1}(J')=G^{i_3}(J)$.
Hence, by (\ref{twof})
$$|G^{i_1}(J)|\le \frac{1}{2}|G^{i_2}(J)|.$$
Using this it follows that the total length of the interval
$J,\dots,G^{k-1}(J)$ contained in one interval $G^{i_3}(J)$
is at most $\sum_{i\ge 0} 2^{-i}=2$ times the length of $G^{i_3}(J)$.
This implies (\ref{satel}). Now (\ref{satel}) gives that
$$\sum_{i=0}^{j-1}|f^i(J)|\le 2n_0$$
where $G^k=f^k$. So if $T\supset J$ is the interval
which is mapped by $f^j$ onto a $\tau/2$-scaled
neighbourhood of $f^k(J)$,
then
$$\sum_{i=0}^{j-1}|f^i(J)|\le S'$$
for some universal constant $S'$.
Here have used the second part of the Koebe Principle.
Thus we have proved (\ref{suff}).
\qed

\subsection{Two step bounds}
For simplicity define
$$d_n=c_{S_n}$$
We shall use boldface letters to indicate the distance to the critical point
(or value), so
  $$\dis{d}{_n} = |d_n-c|, \text{ and } \dis{d}{_n^f}=|d_n^f-c^f|.$$
This notation will also be used for the points we defined before,
namely
$t_n^f$ is the critical point of the monotone
branch of $\fmSn{n}$ near $c^f$ lying on the other side of $c^f$
than $c$ (and therefore than $\cmSn{n}^f$ as well). The critical value
corresponding to $t_n^f$ is $\cSn{n-4}=\fmSn{n}(t_n^f)$.
$$z_n=c_{-S_n}\text{ and }z_n^f=f(z_n)
$$
where $z_n$ could be either to the left or
the right of $c$ depending on the context.
Moreover, remember that we defined
$$y_n=f^{S_n}(c_{S_{n+2}})\text{ and }
y_n^f=f(y_n)$$
in Proposition \ref{top2}.
In the next lemmas the constant $\O$ from
Proposition \ref{crrat} will be written
as $\O_n$, in order to indicate its dependence on $S_n$.
Notice that $\O_n \to 1$ if $n \to \infty$ because of Lemma
\ref{contrpr}.
\begin{lemma} (See \cite{fibo}) \label{lambda}
Let $\lambda^f_n=\dis{d}{^f_{n-2}}/\dis{d}{^f_n}$ then
$\lambda^f_n>3.85$ and $\ln(\dis{d}{^f_{n-4}}/\dis{d}{^f_n})>2.7$
for sufficiently large $n$.
\end{lemma}
\pr
Applying the cross-ratio inequalities
we have
\beqas
\frac{\dis{d}{_n}-\dis{y}{_n}}{\dis{d}{^f_{n+2}}}
\frac{\dis{d}{_{n-4}}}{|t|}
&=&
\frac{|J|}{|j|}\frac{|T|}{|t|}\ge \O_n
\frac{|L|}{|l|}\frac{|R|}{|r|}
\\
&\ge &
\O_n \frac{\dis{y}{_n}}{|\dis{z}{_n^f}-\dis{d}{^f_{n+2}}|}
\frac{\dis{d}{_{n-4}}-\dis{d}{_n}}{|r|}
\eeqas
where $t,j,l,r$ are chosen as in the figure below.
\bigskip

\hbox to \hsize{\hss\unitlength=6mm
\beginpic(20,8)(0,0) \let\ts\textstyle
\put(8,1){{\it Figure 3.2.}}
\put(11,7){\line(1,0){9}}
\put(13,6.8){\line(0,1){0.4}} \put(12.8,6.3){$z_n^f$}
\put(15,6.8){\line(0,1){0.4}} \put(14.8,6.3){$d_{n+2}^f$}
\put(17,6.8){\line(0,1){0.4}} \put(16.8,6.3){$c^f$}
\put(19,6.8){\line(0,1){0.4}} \put(18.8,6.3){$t_n^f$}
\put(16.98,6.8){\line(0,1){0.4}}
\put(12.98,6.8){\line(0,1){0.4}}
\put(14.98,6.8){\line(0,1){0.4}}
\put(18.98,6.8){\line(0,1){0.4}}
\put(13.7,7.3){$l$} \put(16,7.3){$j$} \put(18,7.3){$r$}

\put(0,3){\line(1,0){17}}
\put(1,2.8){\line(0,1){0.4}} \put(0.8,2.3){$c$}
\put(6,2.8){\line(0,1){0.4}} \put(5.8,2.3){$y_n$}
\put(10,2.8){\line(0,1){0.4}} \put(9.8,2.3){$d_n$}
\put(13,2.8){\line(0,1){0.4}} \put(12.8,2.3){$d_{n-4}$}
\put(0.958,2.8){\line(0,1){0.4}}
\put(5.958,2.8){\line(0,1){0.4}}
\put(9.958,2.8){\line(0,1){0.4}}
\put(12.958,2.8){\line(0,1){0.4}}
\put(3.5,3.3){$L$} \put(8,3.3){$J$} \put(11.5,3.3){$R$}

\put(11,6){\vector(-1,-1){2}} \put(8,5){$f^{S_n-1}$}
\endpic\hss}
\bigskip
Using the non-flatness of $c$ and the previous inequality
(and $|t| > |r|$) we get
\beqas
\el \left( 1- \frac{\dis{d}{^f_n}}{\dis{d}{^f_{n-4}}} \right)
&\le&\frac{\dis{d}{_{n-4}}-\dis{d}{_n}}{\dis{d}{_{n-4}}} <
\O_n \frac{\dis{d}{_n}-\dis{y}{_n}}{\dis{y}{_n}}
\frac{|\dis{z}{_n^f}-\dis{d}{^f_{n+2}}|}{\dis{d}{^f_{n+2}}}\cdot 1
\\
&<& \O_n \left(\frac{\dis{d}{_n}}{\dis{d}{_{n+1}}}-1\right)\cdot
\left(\frac{\dis{d}{^f_{n+1}}}{\dis{d}{^f_{n+2}}}-1\right)
\\
&<& \O_n \el \left(\frac{\dis{d}{^f_n}}{\dis{d}{^f_{n+1}}}-1\right)\cdot
\left(\frac{\dis{d}{^f_{n+1}}}{\dis{d}{^f_{n+2}}}-1\right)
\\
&\le&
\O_n \el\left(\sqrt{ \frac{\dis{d}{^f_n}}{\dis{d}{^f_{n+2}}} }-1\right)^2\ .
\eeqas
Here we have used $0<\ell a^{\ell-1}<(b^\ell-a^\ell)/(b-a)<\ell b^{\ell-1}$
and $(\sqrt {ab}-1)^2\ge (a-1)(b-1)$.
Finally we get
\[
1-(\lambda_{n}^f\lambda_{n-2}^f)^{-1}
\le \O_n (\sqrt{\lambda_{n+2}^f }-1)^2
\]
which yields the analogous inequality for $\lambda^f_\infty=
\liminf\lambda^f_n$, and thus $\lambda^f_\infty>3.85$
and $\liminf \ln(\dis{d}{^f_{n-4}}/\dis{d}{^f_n})=2\ln\lambda_\infty^f>2.7$.
One can obtain better estimates for large $\ell$ using
$\ell(\lambda-1)/\lambda\le\ln\lambda^f\le\ell(\lambda-1)$.
\qed
\begin{lemma}{\label{lnln}}
Let $a \in (z_n,\hat z_n)$, $a^f=f(a)$,
$b=\fSn{n}(a)$ and $b^f=f(b)=\fSn{n}(a^f)$.
Then for $n$ large enough
\[
|D\fSn{n}(a^f)|\le \frac{\dis{b}{^f}}{\dis{a}{^f}}
\ln(\frac{\dis{d}{^f_{n-4}}}{\dis{b}{^f}})
\ln(\frac{\dis{d}{^f_n}}{\dis{b}{^f}})
\left(\frac{\dis{d}{^f_{n-4}}}{\dis{b}{^f}}\right)^\el .
\]
\end{lemma} 
\pr
We use the cross-ratio for $\fmSn{n}$ with $l$ is shrunk to a point
$l=\{a^f\}$ and $j=(a^f,c^f)$ and $r=(c^f,t_n^f)$.

\hbox to \hsize{\hss\unitlength=1.3mm
\beginpic(70,50)(0,0) \let\ts\textstyle
\put(7,5){\line(1,0){60}}
\put(7,25){\line(1,0){60}}
\put(7,45){\line(1,0){60}}
\put(28,-3.5){{\it Figure 3.3.}}
\put(15,6){\line(0,-1){2}} \put(14,2){$c$}
\put(30,6){\line(0,-1){2}} \put(29,2){$b$}
\put(45,6){\line(0,-1){2}} \put(44,2){$d_n$}
\put(60,6){\line(0,-1){2}} \put(59,2){$d_{n-4}$}
\put(29,7){$L$} \put(37,7){$J$} \put(52,7){$R$}
\put(15,26){\line(0,-1){2}} \put(14,22){$z_n^f$}
\put(30,26){\line(0,-1){2}} \put(29,22){$a^f$}
\put(45,26){\line(0,-1){2}} \put(44,22){$c^f$}
\put(60,26){\line(0,-1){2}} \put(59,22){$t_n^f$}
\put(29.5,27){$l$} \put(37,27){$j$} \put(52,27){$r$}
\put(15,20){\vector(0,-1){9}}
\put(30,20){\vector(0,-1){9}}
\put(45,20){\vector(0,-1){9}}
\put(60,20){\vector(0,-1){9}}

\put(15,46){\line(0,-1){2}} \put(14,42){$z_n$}
\put(30,46){\line(0,-1){2}} \put(29,42){$a$}
\put(45,46){\line(0,-1){2}} \put(44,42){$c$}
\put(15,40){\vector(0,-1){9}}
\put(30,40){\vector(0,-1){9}}
\put(45,40){\vector(0,-1){9}}

\put(5,15){$f^{S_n-1}$} \put(7,35){$f$}
\endpic\hss}
\vskip0.7cm

\par
In the cross-ratio inequality we can use $|r|<|l|+|j|+|r|$ and have
\beqas
|D\fSn{n}(a^f)|&=&|Df(b)| \, |D\fmSn{n}(a^f)|\\
&\le& \O_n \ell \frac{\dis{b}{^f}}{\dis{b}{}} \,
\frac{\dis{d}{_n}-\dis{b}{}}{\dis{a}{^f}}
\frac{\dis{d}{_{n-4}}-\dis{b}{}}{\dis{d}{_{n-4}}-\dis{d}{_n}}
\\
&=& \O_n \,\frac{\dis{b}{^f}}{\dis{a}{^f}}\,\ell\cdot
\frac{\dis{d}{_n}-\dis{b}{}}{\dis{d}{_n}}
\,\ell\cdot\frac{\dis{d}{_{n-4}}-\dis{b}{}}{\dis{d}{_{n-4}}}
\frac{1}
{\ell\cdot\frac{\dis{d}{_{n-4}}-\dis{d}{_{n}}}{\dis{d}{_n}}}
\,\frac{\dis{d}{_{n-4}}}{\dis{b}{}}
\\
&\le&{\O_n}\,\frac{\dis{b}{^f}}{\dis{a}{^f}}
\,\frac{ \ln\frac{\dis{d}{^f_n}}{\dis{b}{^f}}\ln\frac{\dis{d}{^f_{n-4}}}{\dis{b}{^f}} }
       { \ln\frac{\dis{d}{^f_{n-4}}}{\dis{d}{^f_n}} }
\left(\frac{\dis{d}{^f_{n-4}}}{\dis{b}{^f}}\right)^\el\ ,
\eeqas
which implies the statement as $\O_n \to 1$
and $\ln\frac{\dis{d}{^f_{n-4}}}{\dis{d}{^f_{n}}}>2$.
Here we have used Lemma \ref{lambda}
and the obvious inequalities $\frac{x-1}{x} \le \ln(x)\le x-1$.
\qed
\begin{remark}
We shall use this lemma  several times. In order to simplify the 
notation let us introduce 
$\rho^f_{n}:=\max\{\dis{d}{^f_k}/\dis{d}{^f_{k+1}}\st n-N_0\le k < n\}$,
where $N_0\le 10$ may change from one lemma to another.
Suppose that $b$ in the previous lemma satisfies $\dis{d}{_{n+i}}\le \dis{b}{}$
for some $i\le 6$, then $\dis{d}{^f_n}/\dis{b}{^f}\le (\rho^f_{n+i})^i$,
$\dis{d}{^f_{n-4}}/\dis{b}{^f}\le (\rho^f_{n+i})^{i+4}$  and
\[
|D\fSn{n}(a^f)|\le
\frac{\dis{b}{^f}}{\dis{a}{^f}}\cdot(4+i)\cdot\ln(\rho_{n+i}^f)
\cdot i\cdot\ln(\rho_{n+i}^f)(\rho_{n+i}^f)^{\frac{4+i} {\ell}}\ ,
\]
where in fact $\rho^f$ could have been taken as
$\max\{\dis{d}{^f_k}/\dis{d}{^f_{k+1}}\st n-i-4\le k < n\}$.
\end{remark}
\begin{lemma}\label{derxi}
We have the following estimate
\[
|D\fSn{m}(d_{m+1}^f)| \le 160\cdot\frac{\dis{d}{^f_{m+2}}}{\dis{d}{^f_{m+1}}}
\ln^4(\rho_{m+2}^f)\cdot(\rho_{m+2}^f)^{\frac{13}{\ell}}\ .
\]
\end{lemma}
\pr
We decompose
$D\fSn{m}(d_{m+1}^f)=D\fSn{m-2}(y_{m-1}^f)D\fSn{m-1}(d^f_{m+1})$
and use previous lemma and remark to both factors.
First we put in the lemma  $n=m-2$, $a=y_{m-1}$,
$b=\cSn{m+2}$
and $i=4$ in the remark.
\[
|D\fSn{m-2}(y_{m-1}^f)|\le\frac{\dis{d}{^f_{m+2}}}{\dis{y}{^f_{m-1}}}\cdot 32\cdot
\ln^2(\rho_{m+2}^f)\cdot(\rho_{m+2}^f)^{\frac{8}{\ell}}
\]
Then we put $n=m-1$, $a=d_{m+1}$, $b=y_{m-1}$
and as $d_m \in (y_{m-1},\hat y_{m-1})$ we have $i=1$.
\[
|D\fSn{m-1}(d_{m+1})| \le \frac{\dis{y}{^f_{m-1}}}
{\dis{d}{^f_{m+1}}}\cdot 5\cdot
\ln^2(\rho_{m}^f)\cdot(\rho_{m}^f)^{\frac {5}{\ell}}\ .
\]
The result follows taking  $\rho^f$  depending on  9 consecutive $k$:
$m-7\le k<m+2 $.
\qed
\par
The next lemma prepares the last tool in this subsection.
It describes the estimation (both ways) of $D\fSn{m}(c^f)$.
\begin{lemma}{\label{cone}}
We have for large $m$
\[
\frac{\dis{d}{^f_m}}{\dis{d}{^f_{m+1}}}(\rho_m^f)^{-\frac{4}{\ell}}
\le |D\fSn{m}(c^f)| \le
2\cdot\frac{\dis{d}{^f_m}}{\dis{d}{^f_{m+2}}}\cdot
\ln(\rho_{m+1}^f)\cdot(\rho_{m+1}^f)^{\frac{1}{\ell}}
\]
\end{lemma}
\pr
We use the same trick as in Lemma \ref{lnln}.
\[
|D\fSn{m}(c^f)|=\ell\,\frac{\dis{d}{^f_m}}{\dis{d}{_m}}\,|D\fmSn{m}(c^f)|\ .
\]
For one side we use the cross-ratio for $\fmSn{m}$
on $l=(z_m^f,c^f)$, $r=(c^f,t^f_m)$.

\vskip0.2cm
\hbox to \hsize{\hss\unitlength=1.3mm
\beginpic(60,30)(0,0) \let\ts\textstyle
\put(7,5){\line(1,0){50}}
\put(7,25){\line(1,0){50}}

\put(28,-3){{\it Figure 3.4.}}

\put(15,6){\line(0,-1){2}} \put(14,2){$c$}
\put(30,6){\line(0,-1){2}} \put(29,2){$d_m$}
\put(45,6){\line(0,-1){2}} \put(44,2){$d_{m-4}$}
\put(22,7){$L$} \put(29,7){$J$} \put(37,7){$R$}

\put(15,26){\line(0,-1){2}} \put(14,22){$z_m^f$}
\put(30,26){\line(0,-1){2}} \put(29,22){$c^f$}
\put(45,26){\line(0,-1){2}} \put(44,22){$t_m^f$}
\put(22,27){$l$} \put(29.5,27){$j$} \put(37,27){$r$}
\put(15,20){\vector(0,-1){10}}
\put(30,20){\vector(0,-1){10}}
\put(45,20){\vector(0,-1){10}}

\put(5,15){$f^{S_m-1}$}

\endpic\hss}
\vskip0.5cm
Then we obtain
\beqas
D\fSn{m}(c^f)&\ge&{\O_m}\,\frac{\dis{d}{^f_m}}{\dis{d}{_m}}
\,\ell\,\frac{\dis{d}{_{m-4}}-\dis{d}{_m}}{\dis{d}{_{m-4}}}\,
\frac{\dis{d}{_m}}{\dis{z}{_m^f}}\\
&\ge& \O_m \ln\frac{\dis{d}{^f_{m-4}}}{\dis{d}{^f_m}}\,
\frac{\dis{d}{_m}}{\dis{d}{_{m-4}}}\, \frac{\dis{d}{^f_m}}{\dis{z}{_m^f}}\\
&>& \O_m \frac{\dis{d}{^f_m}}{\dis{d}{^f_{m+1}}}\,
\left(\frac{\dis{d}{^f_m}}{\dis{d}{^f_{m-4}}}\right)^\el\ .
\eeqas
For the other side we take $l=(z_m^f,d_{m+2}^f)$,
$j=(d_{m+2}^f,c^f)$ and $r=(c^f)$.

\hbox to \hsize{\hss\unitlength=1.3mm
\beginpic(60,30)(0,0) \let\ts\textstyle
\put(7,5){\line(1,0){50}}
\put(7,25){\line(1,0){50}}

\put(28,-3){{\it Figure 3.5.}}

\put(15,6){\line(0,-1){2}} \put(14,2){$c$}
\put(30,6){\line(0,-1){2}} \put(29,2){$y_m$}
\put(45,6){\line(0,-1){2}} \put(44,2){$d_m$}
\put(22,7){$L$} \put(37,7){$J$} \put(44,7){$R$}

\put(15,26){\line(0,-1){2}} \put(14,22){$z_m^f$}
\put(30,26){\line(0,-1){2}} \put(29,22){$d_{m+2}^f$}
\put(45,26){\line(0,-1){2}} \put(44,22){$c^f$}
\put(22,27){$l$} \put(37,27){$j$} \put(44,27){$r$}
\put(15,20){\vector(0,-1){10}}
\put(30,20){\vector(0,-1){10}}
\put(45,20){\vector(0,-1){10}}
\put(5,15){$f^{S_m-1}$}
\endpic\hss}
\vskip0.6cm
We obtain, using $\dis{d}{_{m+1}}\le \dis{y}{_m}$,
\beqas
|D\fSn{m}(c^f)| &\le& \O_m \frac{\dis{d}{^f_m}}{\dis{d}{_m}}
\ell\frac{\dis{d}{_m}-\dis{y}{_m}}{\dis{d}{^f_{m+2}}}
\frac{\dis{d}{_m}}{\dis{z}{_m^f}}
\frac{\dis{z}{_m^f}-\dis{d}{^f_{m+2}}}{\dis{y}{_m}}
\\
&\le& \O_m \frac{\dis{d}{^f_m}}{\dis{d}{^f_{m+2}}}
\ell \frac{\dis{d}{_m}-\dis{d}{_{m+1}}}{\dis{d}{_m}}
\frac{\dis{d}{_m}}{\dis{d}{_{m+1}}}
\\
&\le& \O_m \frac{\dis{d}{^f_m}}{\dis{d}{^f_{m+2}}}\ln(\rho_{m+1}^f)
(\rho_{m+1}^f)^{\frac{1}{\ell}}
\ .
\eeqas
And again $\rho^f$ could have been taken as 
$\max\{\dis{d}{^f_k}/\dis{d}{^f_{k+1}}\st m-4\le k\le m \} $.
\qed
\par
\subsection{The  1-step bounds}
We can get a better upper estimate if we combine the previous calculations.  
\begin{prop}\label{onestepbounds}
For $n$ and $\ell$ large enough the derivatives $D\fSn{n}(c^f)$ and the
proportions $\dis{d}{^f_n}/\dis{d}{^f_{n+1}}$ are bounded from above
and separated from below
from 1 by constants independent of $n$ and $\ell$.
\end{prop}
\pr
Consider the following decomposition:
\[
D\fSn{n}(c^f)=D\fSn{n-2}(d_{n-1}^f)\cdot D\fSn{n-1}(c^f)=
D\fSn{n-2}(d_{n-1}^f)D\fSn{n-3}(d_{n-2}^f)\cdot D\fSn{n-2}(c^f).
\]
By Lemma \ref{derxi} used twice with $m=n-2$ and $m=n-3$
in the two first factors
and by Lemma  \ref{cone} used for $m=n-2$ in the third one we have
\[
|D\fSn{n}(c^f)| \le 160^2\cdot 2\cdot\frac{\dis{d}{^f_n}}{\dis{d}{^f_{n-1}}}
\frac{\dis{d}{^f_{n-1}}}{\dis{d}
{^f_{n-2}}}\frac{\dis{d}{^f_{n-2}}}{\dis{d}{^f_n}}\cdot
\ln^{(4+4+1)}(\rho_n^f)(\rho_n^f)^{\frac{13+13+1}{\ell}}
=51200\cdot\ln^9(\rho_n^f)\cdot(\rho_n^f)^{\frac{27}{\ell}}.
\]
This and the other part of Lemma \ref{cone} gives
\[
\frac{\dis{d}{^f_n}}{\dis{d}{^f_{n+1}}}\le
51200\cdot\ln^9(\rho_n^f)\cdot(\rho_n^f)^{\frac{31}{\ell}}\ ,
\]
where again $\rho_n^f$ depends on at most 10 consecutive
quotients $\dis{d}{^f_k}/\dis{d}{^f_{k+1}}$,
with $n-10\le k<n$.
This gives an upper bound of the growth of 
$\rho_\infty^f=\limsup\rho_n^f=\limsup \dis{d}{^f_n}/\dis{d}{^f_{n+1}}$
by
\[
\rho_\infty^f\le
51200\cdot\ln^9(\rho_\infty^f)\cdot(\rho_\infty^f)^{\frac{31}{\ell}}\ ,
\]
(i.e. $\rho^f_\infty<10^{21}$ for large $\ell$)
and proves the upper bound part of the proposition.
For the lower part one can use the estimates from Lemma \ref{lambda},
as from its proof it follows that
$$
\frac{d_n^f}{d_{n+1}^f}  \ge
1 + \frac{1-e^{-2.7}}{\rho_{\infty}^f-1}.
$$
\qed

\begin{prop}\label{kor}
There exists $K > 0$ independent of
$\ell$ and $n$ such that for $\ell$ and $n$
large enough
$$\dis{d}{^f_n}/\dis{u}{_n^f} > 1+K.
$$

\end{prop}
\pr
We shall apply Proposition~\ref{crrat} together with Lemma~\ref{contrpr}
and the remark
after the proof of Lemma \ref{doubc}. Consider $f^{S_{n-1}}$ and
its interval of
monotonicity $t=(z_{n-1}^f,t_{n-1}^f)$ around $c^f$. Let
$l=(z_{n-1}^f,u_{n}^f)$, $j=(u_{n}^f,c^f)$ and $r=(c^f,t_{n-1}^f)$. Denote
by $T,L,J,R$ the images of $t,l,j,r$ under $f^{S_{n-1}}$. Then
\beqas
\frac{\dis{d}{_{n}^f}-\dis{u}{_{n}^f}}{\dis{u}{_{n}^f}}
&\ge&
\frac{\dis{z}{_{n-1}^f}-\dis{u}{_{n}^f}}{\dis{u}{_{n}^f}}=
\frac{|l|}{|j|}\ge
\O_n\frac{|L|}{|J|}\frac{|R|}{|T|}=
\O_n\frac{\dis{u}{_{n-1}^f}}{\dis{d}{_{n-1}^f}-\dis{u}{_{n-1}^f}}
\frac{\dis{d}{_{n-5}^f}-\dis{d}{_{n-1}^f}}{\dis{d}{_{n-5}^f}}\\
&\ge&
\O_n\frac{\dis{d}{_{n}^f}}{\dis{d}{_{n-1}^f}}
\left(1-\frac{\dis{d}{_{n-1}^f}}{\dis{d}{_{n-5}^f}}\right)
\ge
\O_n\frac{1-\frac{1}{2.7}}{\rho_n^f}
\ ,
\eeqas
and this is bounded away from $0$ uniformly in $n$.
\qed

The main result of this section is the finally the following
theorem.

\begin{theo}\label{realb}
$$
\frac{\dis{d}{_n^f}}{\dis{u}{_n^f}}\,\,,\,\,
\frac{\dis{d}{_n^f}}{\dis{d}{_{n+1}^f}}\text{ and }
\frac{\dis{u}{_n^f}}{\dis{u}{_{n+1}^f}}$$
are bounded and bounded away from one for all $\ell$
and $n$ large enough. In particular, there are constants
$C_1,C_2$ for which
$$
\frac{C_1}{\ell}
\le \,\,
 \frac{|d_n-u_n|}{|u_n-c|}\,\,,\,\,
\frac{|\dis{d}{_n}-\dis{d}{_{n+1}}|}{|d_n-c|}\,\,,\,\,
\frac{|\dis{u}{_n}-\dis{u}{_{n+1}}|}{|u_n-c|} \,\,
\le \frac{C_2}{\ell}.
$$
\end{theo}
\pr Follows from the previous two results and the fact that
$f$ has a critical point of order $\ell$ at $c$. \qed

\vskip 0.5cm
\hbox to \hsize{\hss\unitlength=5mm
\beginpic(20,2)(0,0) \let\ts\textstyle
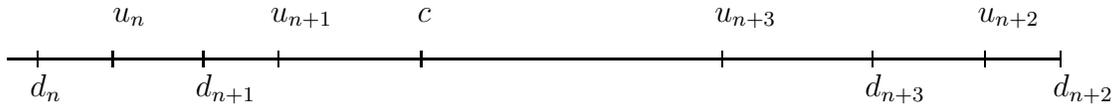
\put(-1,-2.5){{\it Figure 3.6: The points $u_n$ and $d_n$.
are on the same side of $c$. The points}}
\put(-1,-4){{\it $d_{n+2}$ and $d_n$
are on opposite sides of $c$;
$z_{n-1}$ is between $d_n$ and $u_n$.}}
\put(-4,0.2){\line(1,0){28}}
\put(-3.2,0){\line(0,1){0.4}} \put(-3.4,-0.8){$d_n$}
\put(-1.2,0){\line(0,1){0.4}} \put(-1.2,1.2){$u_n$}
\put(1.2,0){\line(0,1){0.4}} \put(1,-0.8){$d_{n+1}$}
\put(3.2,0){\line(0,1){0.4}} \put(3,1.2){$u_{n+1}$}
\put(7,0){\line(0,1){0.4}} \put(6.9,1.2){$c$}
\put(15,0){\line(0,1){0.4}} \put(14.8,1.2){$u_{n+3}$}
\put(19,0){\line(0,1){0.4}} \put(18.8,-0.8){$d_{n+3}$}
\put(22,0){\line(0,1){0.4}} \put(21.8,1.2){$u_{n+2}$}
\put(24,0){\line(0,1){0.4}} \put(23.8,-0.8){$d_{n+2}$}
\endpic\hss}
\vskip 3cm

\sect{The random walk argument}
\newcommand{\jmax}{d}
\newcommand{\kmax}{k_0}
\newcommand{\numax}{\nu_{{\rm max}}}
In this section we shall state and prove an abstract result
about the evolution of typical points
under a (nearly) Markov map with a kind of random walk structure.
So
let $(X,\F,m)$ be some space with probability measure $m$
and $\sigma$-algebra $\F$. Let $\A=\{A_k\colon k=0,1,2,\ldots\}$ 
denote a partition of
$X$ into $\F$-measurable sets, and let
$F\colon X\to X$ be a $\F$-measurable transformation. 
We denote $\A_n=\bigvee_{k=0}^{n-1}F^{-k}\A$ and
let $\H$ be
the family of all measures of the form 
$F_*^n(m_{|A})$ with $A\in\A_{n+1}$ and $n\ge 0$.
Now take $r,\kmax\in \nz$ and define  for $i\ge 1$,
\beq
\label{deff}
a_i=m(A_{r+i-\kmax-1})\mbox{ and }
\nu_i=\mu((F_{|A_r})^{-1}(A_{r+i-\kmax-1}))/\mu(A_r)\ ,
\eeq
where $\mu$ is some measure from the class $\H$ defined above.
Note that
$$
\frac{a_{k+1}}{a_{j+1}}
\frac{\nu_{j+1}}{\nu_{k+1}}$$
is equal to one if $F$ preserves the measure $\mu$.
 
\par
In the remainder of this section the sequences of
positive real numbers $(a_i)_{i\ge0}$ and
$(\nu_i)_{i\ge 1}$ will be assumed to have
a particular exponential decay.
Here $a_i$, $\nu_i$ are as above
and $a_0$ will be a suitable constant corresponding to
a constant which comes from `Koebe space'.
So we say that two such sequences satisfy
the {\em scaling
condition}
with constants $\rhom,\rhop,\O_1,\O_2>0$ and $\jmax\ge 0$, if 
\beqa
\label{tfirst}
\rhom^{k-j}&\le&\frac{\sumai{j}}{\sumai{k}}\le\rhop^{k-j}
\quad\mbox{ for all }0\le j\le k\ ,\\
\label{tsecond}
\O_1 \cdot \frac{a_0}{\sumaI{k}}\cdot\frac{\numax}{\sumaI{k+1}}
&\le&
\frac{\nu_{k+1}}{a_{k+1}}
\quad\mbox{for all $k\ge d$,}\\ 
\label{tthird}
\frac{\nu_{k+1}}{a_{k+1}}
&\le&
\O_2 \cdot \frac{\numax}{a_1}\quad\mbox{for all $k\ge 0$,} 
\eeqa
where $\numax=\max\{\nu_1,\ldots,\nu_{\jmax+1}\}$.
As we shall see in the lemma below, (\ref{tfirst}) means
that the numbers $a_i$ decay at a slow rate if $\rho_i>1$
are close to one. Moreover, this lemma implies that
if $(\rho_1-1)/(\rho_0-1)$ is not too large then
(\ref{tsecond})
is equivalent to the more symmetric expression 
$$
\O_1 \cdot
\frac{\sumaI{d}}{\sumaI{k}}\cdot\frac{\sumaI{d+1}}{\sumaI{k+1}}
\cdot \frac{a_{k+1}}{a_{\max}}
\le 
\frac{\nu_{k+1}}{\numax}.
$$
where $a_{\max}=\max\{a_1,\dots,a_{d+1}\}$.
The previous inequalities
(\ref{tfirst}) and (\ref{tsecond}) combined show that
the ratio of the `mass'
going from state $A_r$ to state $A_{r+i-k_0-1}$
compared to the mass going to one of the states $A_{r-k_0},
\dots,A_{r-k_0+d}$
goes only down slowly with $i$. So - roughly speaking -
a reasonably large set of points move to a state with much larger
index. This suggests that $m$-typical points will move to states
with larger and larger indices.
This intuitive idea is formalized in the following theorem.
\par
\begin{theo}\label{random}
Suppose there is $\kmax\in\nz$ such 
that $F(A_r)\subseteq\cup_{j=0}^{\infty}A_{r-\kmax+j}$ for all
$r\ge 2$.
Then there exists for each $C>1$, $\O_1,\O_2>0$, and $\jmax\in\nz$
a constant $\rho\in(1,2)$ with the property that if the
assumption 
stated below is satisfied, 
then there exists a set $D\in\F$ with $m(D)>0$
such that
for each $x\in D$ and each $A_j$ the set
\[
\{k: F^k(x)\in A_j\}
\]
has finite cardinality. The assumption is: 
\begin{quote}
For any $1<\rhom<\rhop<\rho$ with
$\frac{(\rhop-1)}{(\rhom-1)}\le C$
there is $r_0>0$ such that
for any $r\ge r_0$ and any $\mu\in\H$ with $\mu(A_r)>0$
there exists $a_0>0$ such that the sequences
$(a_i)_{i\ge 0}$ and $(\nu_i)_{i\ge 1}$ satisfy the {\em scaling
condition}
with constants $\rhom,\rhop,\O_1,\O_2$ and $\jmax$, 
where $a_i$ and $\nu_i$ ($i\ge 1$) are defined as in (\ref{deff}).
\end{quote}
\end{theo}
Before we turn to the proof of the theorem in the next three
subsections, we
state some simple properties of sequences satisfying the 
{\em scaling condition}. These properties give a better intuition
for the meaning  of this condition.
\begin{lemma}\label{lprop}
Let $(a_i)_{i\ge 0}$, $(\nu_i)_{i\ge 1}$ satisfy condition
(\ref{tfirst}) and define
$\Kp=(\rhop-1)\rhom/(\rhom-1)\rhop$.
Then for any $0\le j\le k$ we have
\beqa
1-\rhom^{-1}&\le&\frac{a_j}{\sumai{j}}\le 1-\rhop^{-1}\label{li}
\\
\frac{1}{\Kp}\rhom^{k-j}&\le&\frac{a_j}{a_k}\le\Kp\rhop^{k-j}\
.   \label{lii}
\eeqa
If also condition (\ref{tthird}) is satisfied, then for $k\ge 1$,
\beq
\nu_k\le \O_2\numax \Kp\rhom^{-(k-1)}\ .
\label{lvii}
\eeq
\end{lemma}
\proof
The proof follows immediatly by calculation. For example
(\ref{li}):
\[
\frac{a_j}{\sumai{j}}=\frac{\sumai{j}-\sumai{j+1}}{\sumai{j}}\le
1-\rhop^{-1}\ .
\]
\qed
\subsection{The martingale argument}
%
As before 
let $(X,\F,m)$ be some space with probability measure $m$
and $\sigma$-algebra $\F$ and $\A=\{A_k\colon k=0,1,2,\ldots\}$ 
a partition of
$X$ into $\F$-measurable sets. 
$F\colon X\to X$ is a $\F$-measurable transformation, 
$\A_n=\bigvee_{k=0}^{n-1}F^{-k}\A$, and 
$\H=\{F_*^n(m_{|A})\colon A\in\A_n, n\ge 0\}$.
\par
Observe that $\A$ is a Markov 
partition for $F$ if and only if $F^kA$ is an element of 
$\A$ for each $A\in\A_{k+1}$ and each $k\ge 0$.
In order to make the following 
proposition most widely applicable we shall not assume that $F$
is strictly 
Markov but formulate instead some restrictions on $\H$.
Furthermore, even if $F$ is topologically 
Markov, the nonlinearity of its branches still prevents $F$ to
be also 
measure theoretically Markov. Therefore we do not use in our
proof 
a Markov-like 
model but instead a more flexible martingale construction. 
As a general reference to the theory of martingales we give 
\cite{stout}.
\par
Define $\ph:X\to\{0,1,2,\ldots\}$ by
\[
\ph(x)=n\mbox{ if }x\in A_n
\]
and
\[
\Delta\ph:=\ph\circ F-\ph\ .
\]
\begin{prop}\label{martin}
Assume there are $r_0\in\nz$ and $M>0$ such that for
any $A\in\A_{k+1}$, $k\ge 0$, with $\ph_{|F^kA}\ge r_0$ holds:
\beqa\label{drift1}
\int_A(\Delta\ph-1)\circ F^k\, dm &\ge& 0\quad\mbox{ and}
\\ 
\label{secmoments}
\int_A(\Delta\ph)^2\circ F^k\, dm &\le& M\cdot m(A)
\quad\mbox{ for all } n\ge0\ .
\eeqa
Then
\[
\liminfn\frac{\ph\circ F^n}{n} \ge 1
\quad\mbox{$m$-a.s.,}
\]
and there exists a set $D\in\F$ with $m(D)>0$
such that for every $x\in D$ the trajectory $x,Fx,F^2x,\ldots$
visits
each set $A_k\in\A$ only finitely often.
\end{prop}
\pr
Fix $s>r_0$ and denote by $\mu$ the normalized restriction of $m$
to 
$A_s$.
Let $\F_n$ be the $\sigma$-algebra
generated by the partition $\A_{n+1}$.
Then $\ph \circ F^n$ is $\F_n$-measurable, i.e.,
\[
E_{\mu}[\ph\circ F^n|\F_n]=\ph\circ F^n\ .
\]
Define a stopping time $\tau:X\to\nz\cup \{\infty\}$ by
\[
\tau(x)=\alternative
{\infty}{\mbox{ if }\ph(F^nx)>r_0\mbox{ for all }n\ge 0}
{\min\{n\ge 0\colon \ph(F^nx)\le r_0\}}{\mbox{ otherwise,}}
\]
and the random variables $(Z_n)_{n\ge 0}$ by 
\[
Z_n(x)=\alternative
{\ph(F^nx)}{\mbox{ if }\tau(x)>n}
{\ph(F^{\tau(x)}x)}{\mbox{ if }\tau\le n.}
\]
Then also the $Z_n$ are $\F_n$-measurable.
So for any $x\in X$ and $n\ge 0$ with $\tau(x)>n$,
\beqa
&&E_{\mu}[Z_{n+1}|\F_n](x)-Z_{n}(x)-1\nonumber\\
&=&
E_{\mu}[(\Delta\ph-1)\circ F^n|\F_n](x)
\label{emu}\\ 
&=&
\sum_{A\in\A_n}\chi_A(x)\cdot\frac
1{m(A)}\int_A(\Delta\ph-1)\circ F^n\, d\mu
\nonumber\\
&\ge& 0\ ,
\nonumber
\eeqa
where we used  
(\ref{drift1}) for the inequality. If $\tau(x)\le n$, then 
$E_{\mu}[Z_{n+1}|\F_n](x)-Z_{n}(x)=0$. Note that in both cases
$E_{\mu}[Z_{n+1}|\F_n](x)\ge Z_{n}(x)$, i.e. $(Z_n,\F_n)_{n\ge
0}$ is a 
submartingale with respect to $\mu$.
Now define
\[
W_n=Z_0+\sum_{k=1}^n(E_{\mu}[Z_k|\F_{k-1}]-Z_{k-1})
\quad\mbox{ and }\quad M_n=Z_n-W_n
\]
(this is, by the way, the Doob-decomposition of $(Z_n,\F_n)_{n\ge
0}$).
Then $W_0=Z_0=s$ and $M_0=0$ $\mu$-a.s., and $(M_n,\F_n)_{n\ge
0}$ is a martingale:
\[
M_{n+1}-M_n=Z_{n+1}-Z_n-W_{n+1}+W_n=
Z_{n+1}-E_{\mu}[Z_{n+1}|\F_n]
\]
and therefore
\[
E_{\mu}[M_{n+1}|\F_n]=
E_{\mu}[M_n|\F_n]+E_{\mu}[Z_{n+1}|\F_n]-E_{\mu}[E_{\mu}[Z_{n+1
}|\F_n]|\F_n]
=M_n\ .
\]
$(W_n,\F_{n-1})_{n\ge 1}$ is a predictable stochastic
sequence with
\[
W_{n+1}-W_n=E_{\mu}[Z_{n+1}|\F_n]-Z_n\ge \chi_{\{\tau>n\}}
\]
because of (\ref{emu}). It follows that 
\beq\label{W_ngrowth}
W_n\ge n+s\quad\mbox{ on }\{x\st \tau(x)\ge n\}\ .
\eeq 
Next note that on $\{\tau>n\}$ holds
\beqas
&&E_{\mu}[(M_{n+1}-M_n)^2|\F_n]
=
E_{\mu}[(Z_{n+1}-E_{\mu}[Z_{n+1}|\F_n])^2|\F_n]\\
&\le&
E_{\mu}[(Z_{n+1}-Z_{n})^2|\F_n]
=
E_{\mu}[(\Delta\ph)^2\circ F^n|\F_n]\\
&\le&
M\ ,
\eeqas
where we used the fact that $E_{\mu}[(Z_{n+1}-Y)^2|\F_n]$ 
is minimized by $Y=E_{\mu}[Z_{n+1}|\F_n]$ for the first
inequality and
assumption (\ref{secmoments}) for the second one.
On $\{\tau\le n\}$ we have
\[
M_{n+1}-M_n=Z_{n+1}-E_{\mu}[Z_{n+1}|\F_n]=Z_{n+1}-Z_n=0\ .
\]
Both estimates together yield 
$E_{\mu}[(M_{n+1}-M_n)^2]\le M$, and we can apply Chow's version
of the 
Hajek-R\'enyi inequality (see \cite[Theorem 3.3.7]{stout}):
\beq\label{chow_est}
\mu\left\{\max_{1\le i\le n}\frac{|M_i|}{s-r_0+i}\ge 1 \right\}
\le
\sum_{i=1}^{n}\frac{E_{\mu}[(M_{i}-M_{i-1})^2]}{(s-r_0+i)^2}
\le
M\cdot\sum_{j>s-r_0}\frac 1{j^2}<\frac 12\ ,
\eeq
if $s-r_0$ is large enough. Hence
\beqas
\mu\{\tau<\infty \}&=&
\mu\left(\bigcup_{n\ge 1}\{\tau=n\}\right)
\le
\mu\left(\bigcup_{n\ge 1}\{Z_n\le r_0\mbox{ and } W_n\ge n+s\}
\right)\\
&\le&
\mu\left(\bigcup_{n\ge 1}\{M_n\le r_0-s-n\}\right) 
\le
\mu\left(\bigcup_{n\ge 1}\{|M_n|\ge s-r_0+n\}\right) \\
&=&\sup_{n\ge 1}\mu\left\{\max_{1\le i\le
n}\frac{|M_i|}{s-r_0+i}\ge 1\right\}
\le
\frac 12
\eeqas
for such $s$, i.e.\ $\mu\{\tau=\infty\}\ge\frac 12$.
\par
Now a convergence theorem of Chow (see \cite[Theorem
3.3.1]{stout}) asserts 
that 
\[
\limn M_n/(s-r_0+n)=0\quad\mbox{ $\mu$-a.s.}
\]
in view of the finiteness of the 
sum in (\ref{chow_est}). Hence, on $\{\tau=\infty\}$,
\[
 \liminfn\frac{\ph\circ F^n}{n}=\liminfn\frac{Z_n}{n}
=\liminfn\frac{W_n}{n}+\limn\frac{M_n}{n}\ge 1
\quad\mbox{$\mu$-a.s.}
\]
in view of (\ref{W_ngrowth}).
In particular,
for each $x\in \{\tau=\infty\}$ 
the trajectory $x,Fx,F^2x,\ldots$ visits each element $A_k\in\A$
only
finitely often.
\qed
\subsection{Some calculations}
Let $(a_i)_{i\ge0}$ and
$(\nu_i)_{i\ge 1}$ be two sequences of positive real numbers
which 
satisfy the {\em scaling condition}
(\ref{tfirst},\ref{tsecond},\ref{tthird})
with constants $\rhom,\rhop,\O_1,\O_2,\jmax$. We assume 
additionally that
$\sum_{i=0}^\infty a_i<\infty$ and
$\sum_{i=1}^\infty \nu_i =1$,
\begin{prop}\label{expv}
For any $E, \O_1,\O_2,C,\jmax$ there is a $\rho\in (1,2)$
such that
if the numbers $\rhom, \rhop$ from above satisfy
$1<\rhom<\rhop<\rho$ and
$\frac{(\rhop-1)}{(\rhom-1)}\le C$ 
then
\beq\label{firstmoment}
\sum_{j=1}^{\infty} j\nu_j > E\ ,
\eeq
and there is some constant $M>0$ depending only on $\O_2,C,\rhom$
such 
that
\beq\label{secondmoment}
\sum_{j=1}^{\infty} j^2\nu_j < M\ .
\eeq
\end{prop}
For the proof we need several lemmas.
\begin{lemma}\label{gerlem}
Let $(q(j))$ be a positive increasing sequence.
Then, if $n-1\ge d+1$,
\[
\sum_{j=d+1}^{n-1} j\cdot\frac{q(j+1)-q(j)}{q(j)q(j+1)}
\ge
\sum_{j=d+1}^{n-1} \frac{q(n)-q(j)}{q(n)q(j)}\ .
\]
If $q(\infty):=\limn q(n)$ exists and if 
$\limn n\cdot(q(\infty)-q(n))=0$, then the above inequality holds
also for $n=\infty$.
\end{lemma}
\proof
\beqas
\sum_{j=d+1}^{n-1} j\cdot\frac{q(j+1)-q(j)}{q(j)q(j+1)}
&=&
-\sum_{j=d+1}^{n-1} j\cdot\left(\frac
{1}{q(j+1)}-\frac{1}{q(j)}\right)\\
&=&
-\sum_{j=d+1}^{n-1}\frac j{q(j+1)}+\sum_{j=1}^{n-1}\frac
j{q(j)}\\
&=&
\sum_{j=d+1}^{n-1}\frac 1{q(j)}(j-(j-1))
+\frac{d}{q(d+1)}-\frac {n-1}{q(n)}\\
&=&
\sum_{j=d+1}^{n-1}\left(\frac{1}{q(j)}-\frac{1}{q(n)}\right)
+d\cdot\left(\frac{1}{q(d+1)}-\frac{1}{q(n)}\right)\\
&\ge&
\sum_{j=d+1}^{n-1} \frac{q(n)-q(j)}{q(n)q(j)} \ .
\eeqas
As, under the additional assumption, $\limn q(n)=q(\infty)$ and 
\[
\sum_{j=1}^{n-1} \frac{q(\infty)-q(n)}{q(n)q(j)}
\le
\frac{n}{q(1)^2}(q(\infty)-q(n))
\to 0\quad\mbox{as }n\to\infty\ ,
\]
also the inequality for $n=\infty$ follows.
\qed
\begin{lemma}\label{integr}
Let $1<\rho<2^{1/d}$. Then
\[
\sum_{k=d+1}^{\infty}\frac{1}{\rho^k-1}>
\frac{1}{\ln \rho}\ln\left(\frac{1}{2(d+1)(\rho-1)}\right)\ .
\]
\end{lemma}
\proof
Fix $M>0$ such that $\rho^{M+1}-1>\rho^M$. Then
\beqas
\sum_{k=d+1}^M\frac{1}{\rho^k-1}
&>&
\int_{d+1}^{M+1}\frac{1}{\exp(x\ln\rho)-1}dx\\
&=&
-(M-d)+\frac{1}{\ln\rho}
\left(\ln (\rho^{M+1}-1)-\ln(\rho^{d+1}-1)\right)\\
&\ge&
d+\frac{1}{\ln\rho}\ln\left(\frac{1}{\rho-1}\cdot
\frac{\rho-1}{\rho^{d+1}-1}\right)\\
&\ge&
\frac{1}{\ln\rho}\ln\left(\frac{1}{2(d+1)(\rho-1)}\right).
\\
\eeqas
Here we have used $1<\rho<2^{1/d}$ in the last inequality.
\qed
\begin {lemma}\label{bigmu}
Let $(a_i)$, $(\nu_i)$ and all constants be as in Proposition
\ref{expv}.
Then 
\[
\sum_{k=1}^{\infty} k\,\nu_k\ge
\frac{\numax}{2}\cdot
\O_1\cdot \frac{1}{\Kp}\ln\left(\frac{1}{2(d+1)(\rhop-1)}\right)
\]
provided $1<\rhop<2^{1/d}$,
where $\Kp$ was defined in Lemma \ref{lprop}.
\end{lemma}
\proof
Let $q(j)=\sum_{i=0}^{j-1}a_i$. As the $a_i$ decrease
exponentially by 
(\ref{lii}), $q(\infty)=\lim_{j\to\infty}q(j)$ exists, and 
we can apply Lemma \ref{gerlem}:
\beqas
&&\sum_{k=d+1}^{\infty} k\cdot\frac{q(k+1)-q(k)}{q(k)q(k+1)}
\ge\sum_{k=d+1}^{\infty}\frac{q(\infty)-q(k)}{q(\infty)q(k)}
\\
&=&
\frac{1}{q(\infty)}
\sum_{k=d+1}^{\infty}\frac{\sumai{k}}{\sumai{0}-\sumai{k}}
=
\frac{1}{q(\infty)}
\sum_{k=d+1}^{\infty}\frac{1}
{\sumai{0}/\sumai{k}-1}
\\
&\ge&
\frac{1}{q(\infty)}
\sum_{k=d+1}^{\infty}\frac{1}{\rhop^{k}-1}
\ge
\frac{1}{\sumai{0}}
\frac{1}{\ln \rhop}\ln\left(\frac{1}{2(d+1)(\rhop-1)}\right)\ .
\eeqas
For the last two inequalities we have used (\ref{tfirst})
and Lemma \ref{integr}.
Observe next that by (\ref{tsecond}) we have 
\[
\frac{\nu_k}{\numax}\ge\O_1 \cdot
\frac{a_0}{\sumaI{k-1}}\cdot\frac{a_{k}}{\sumaI{k}} 
= \O_1\cdot
a_0\cdot\frac{q(k+1)-q(k)}{q(k)q(k+1)}\ .
\]
Combining this with the previous estimate we obtain
\beqas
\sum_{k=1}^{\infty} k\,\nu_k 
&=&
\numax\sum_{k=d+1}^{\infty} k\,\frac{\nu_k}{\numax}\\
&\ge&
\numax\cdot \O_1\cdot a_0\cdot\sum_{k=d+1}^{\infty} k\cdot
\frac{q(k+1)-q(k)}{q(k)q(k+1)}
\\
&\ge&
\numax\cdot\O_1 \cdot \frac{a_0}{\sumai{0}}\cdot
\frac{1}{\ln \rhop}\cdot\ln\left(\frac{1}{2(d+1)(\rhop-1)}\right)
\\
&\ge&
\numax\cdot \O_1 \cdot
\frac{\rhom-1}{\rhom\ln
\rhop}\cdot\ln\left(\frac{1}{2(d+1)(\rhop-1)}\right)\\
&\ge&
\frac{\numax}{2}\cdot \O_1\cdot
\frac{1}{\Kp}\frac{\rhop-1}{\ln\rhop}
\ln\left(\frac{1}{2(d+1)(\rhop-1)}\right)
\ .
\eeqas
For the last two inequalities we have used (\ref{li}), $\rhop<\rho<2$
and the definition of $\Kp$ in Lemma \ref{lprop}. As 
$\ln(\rhop)\le\rhop-1$, this proves the lemma.
\qed
\begin{lemma}\label{smallmu}
Let $(a_i)$, $(\nu_i)$ and all constants be as in Proposition
\ref{expv},
and let
$\Kp$
be defined as in Lemma \ref{lprop}.
Then for any $r>0$ we have
\[
\sum_{k=1}^{\infty} k\,\nu_k
>r (1-\O_2\numax\Kp\,r)
\ .
\]
\end{lemma}
\proof
The idea is to use the very rough estimation
\[
\sum_{k=1}^{\infty} k\,\nu_k>\sum_{k=r}^{\infty} k\,\nu_k>
r(1-\sum_{k=1}^{r-1}\nu_k)\ .
\]
Using (\ref{lvii}),
\[
\sum_{k=1}^{r-1}\nu_k\le
\O_2\numax \Kp \sum_{k=1}^{r-1}\rhom^{-(k-1)}
\le\O_2\numax\Kp\,r
\]
the lemma follows immediately.
\qed

\noindent
{\em Proof of Proposition~\ref{expv}:}\quad
Recall that $\Kp=(\rhop-1)\rhom/(\rhom-1)\rhop$.
We have that $\Kp\le C$.
Fix $r=2E$ and choose $\rho\in (1,2)$ such that
\[
\ln\frac{1}{2(d+1)(\rho-1)}=64\frac{\O_2}{\O_1} E^2C\ .
\]
Then 
$\ln\frac{1}{2(d+1)(\rhop-1)}>64\frac{\O_2}{\O_1} E^2C$, and
by Lemma \ref{smallmu} we have
\beqas
\sum_{k=1}^{\infty} k\,\nu_k&>&2E(1-2E\numax \O_2 2C)
\\
&>&
2E\left(1-\frac{\O_1}{16EC}\numax\cdot \ln
(\frac{1}{2(d+1)(\rhop-1)})\right)\ ,
\eeqas
which is bigger than $E$ if
$\numax\ln\frac{1}{2(d+1)(\rhop-1)}<8EC/\O_1$. 
Otherwise $\sum_{k=1}^{\infty} k\,\nu_k>E$ follows from Lemma
\ref{bigmu}.
The existence of a uniform bound (\ref{secondmoment}) 
for the second moments of $(\nu_k)$ 
follows from (\ref{lvii}).
\qed
\subsection{The proof of Theorem 4.1}  
%
In this section we shall prove Theorem \ref{random}.
Take $\mu\in\H$, i.e.\ fix some $A\in\A_{n+1}$ and 
consider the measure $\mu=F_*^n(m_{|A})$.
Define $\nu$ to be the normalization of $\mu$ on $A_r$, i.e.
$\nu(\cdot)=\mu(\cdot \cap A_r)/\mu(A_r)$.
For $j\ge 1$ define
\[
a_j=m(A_{r+j-\kmax-1})\mbox{ and
}\nu_j=\nu((F_{|A_r})^{-1}(A_{r+j-\kmax-1}))\ .
\]
Let $E=\kmax+2$. If the assumptions of Theorem \ref{random} are
satisfied, then
we get from Proposition \ref{expv} that
\[
\sum_{j=1}^{\infty} j\nu_j>\kmax+2\quad\mbox{ and }\quad
\sum_{j=1}^{\infty} j^2\nu_j<M
\]
where $M$ does not depend on the particular measure $\mu$.
Hence
\[
\sum_{j=1}^{\infty} j\nu(F_{|A_r})^{-1}(A_{r+j-\kmax-1}) >\kmax+2\
.
\]
Observe that $\Delta\ph$ from Proposition \ref{martin} 
is equal to $j-\kmax-1$ on
$(F_{|A_r})^{-1}(A_{r+j-\kmax-1})$.
So
\beqas
\int_{A_r}(\Delta\ph-1)\, d\nu&=&
\left[\sum_{j=1}^\infty
(j-\kmax-1)\nu(F_{|A_r})^{-1}(A_{r+j-\kmax-1})\right] -1\\
&\ge&
\left[\sum_{j=1}^\infty j\nu(F_{|A_r})^{-1}(A_{r+j-\kmax-1})\right]
-\kmax-2
\ge 0 \ ,
\eeqas
and similarly 
\[
\int_{A_r}(\Delta\ph)^2\, d\nu<M\ .
\]
Hence the assumptions of Proposition \ref{martin}
are satisfied and this implies that the
assertion of Theorem \ref{random} holds.

\sect{Proof of the Main Theorem}
In this section we shall complete the proof of the Main Theorem.
So let $f$ be a $C^2$ Fibonacci map with
a critical point of order $\ell$.
First we should remark that the complement of the
basin of $\omega(c)$ is a residual set.
This can be seen as follows. From Chapter IV of \cite{MS}
it follows that $f$ has no wandering intervals
(a wandering interval is an interval whose
forward iterates are all disjoint and which is not in the
basin of a periodic attractor).
Moreover, $f$ is not renormalizable
and has positive topological entropy, see \cite{HK}.
It follows that $f$ is semi-conjugate to a tent-map of the form
$$
x\mapsto \lambda\left( 1- |2x-1|\right)
$$
and that the semi-conjugacy only collapses
components of basins of periodic attractors.
Clearly such components cannot be in the basin of the Cantor set
$\omega(c)$. So it suffices to show that there exists
a residual set of points $x$ for which $\omega(x)$ (w.r.t.
a tent-map) is equal to a cycle of intervals. This fact
is well-known, see for example \cite[page 189]{Mil}.
So the deepest part of the proof consists in showing
that $B(\omega(c))$ has positive Lebesgue measure.
\par
Let the points $u_k$, $c_{S_k}$ and so on be defined as in Section
\ref{sectop} and choose as before
$\tilde u_{k+1}\in \{u_{k+1},\hat u_{k+1}\}$
so that it is on the same side of $c$ as $u_k$.
Define intervals
$I_k=(u_k,\tilde u_{k+1})$ and $\hat I_k$ (the interval symmetric
to $I_k$),
and a map
\[
F\,\,\colon \,\,\bigcup (I_k \cup \hat I_k) \to
\bigcup (I_k\cup \hat I_k)
\]
by
\[
F|I_k=f^{S_k}\ .
\]
Then for $k>1$
\[
F(I_k)=F(\hat I_k)=(u_{k-2},u_k)\ .
\]
Hence, if we let $A_k=I_k\cup\hat I_k$ ($k\ge 0$), then 
$\A=\{A_k\colon k=0,1,2,\ldots\}$ is a partition of $X=(u_0,\hat
u_0)$,
and $F$ is Markov with respect to $\A$. 

\hbox to \hsize
{\hss\unitlength=3mm
\beginpic(45,55)(0,0) \let\ts\textstyle
\put(5,30){\line(1,0){30}}
\put(10,2){\line(0,1){50}}
\put(15,0){{\it Figure 5.1.}}
\put(17,30){\squine(0, -0.5, -2, 0, 9, 10)}
\put(17,30){\squine(0, 1, 3, 0, -16, -17)}
\put(23,30){\squine(0, -1, -3, 0, 16, 17)}
\put(23,30){\squine(0, 1, 4, 0, -25.5, -27)}

\put(9,29){$c$}
\put(15,30){\line(0,1){10}} \put(14.2,28.5){$\tilde u_{k+1}$}
\put(20,30){\line(0,1){17}} \put(20,30){\line(0,-1){17}}
\put(20.5,28.5){$u_k$}
\put(27,30){\line(0,-1){27}} \put(27.5,28.5){$u_{k-1}$}
\put(17.5,28.5){$z_k$} \put(23.5,28.5){$z_{k-1}$}

\put(10,3){\line(1,0){17}} \put(7,3){$u_{k-3}$}
\put(10,13){\line(1,0){10}} \put(7,13){$u_{k-2}$}
\put(10,40){\line(1,0){5}} \put(8,40){$u_k$}
\put(10,47){\line(1,0){10}} \put(7,47){$u_{k-1}$}

\put(15,25){$f^{S_k}$} \put(24,37){$f^{S_{k-1}}$}

\endpic\hss}

In this section we shall show that $F$ and $\A$ satisfy the
assumptions
of Theorem \ref{random} with $\jmax=2$, $\kmax=2$ and thus prove
\begin{theo}\label{drift}
For all sufficiently large $\ell$ holds:
The set $D$ of all points $x$
for which the trajectory $(F^kx)_{k>0}$ visits each interval
$I_n$ and $\hat
I_n$ at most finitely often, has positive Lebesgue measure. 
\end{theo}
Let us first show that this result implies our Main Theorem,
which states:
\begin{theo}
$\omega_f(c)$ is an absorbing
Cantor set attractor for $f$ provided $\ell$ is large enough.
\end{theo}
\pr
First we should remark that $f$ is ergodic with respect to
the (non-invariant) Lebesgue measure if its Schwarzian derivative
is
negative, see \cite{BL}.
Hence $F$ is ergodic with respect to Lebesgue measure, and as
$f^{-1}(D)=D$ and $D$ has positive Lebesgue
measure, $D$ has full Lebesgue measure in this case.
If $f$ is a smooth Fibonacci map with a periodic attractor
then, of course, $f$ is not ergodic and $\omega(c)$ cannot be
a `global' attractor. But even in this case, the argument below
will show
that it attracts a set of points of positive Lebesgue measure.
We should remark that $\omega(c)$ is not accumulated by
periodic attractors, see \cite{MMS} or \cite{MS}, so near
$\omega(c)$
these periodic attractors are `invisible'.
\par
So consider a point $x\in X$ for which $(F^kx)_{k>0}$
visits each interval $I_n$ 
and $\hat I_n$ at most finitely often, and denote by
$t_1<t_2<t_3<\ldots$ the 
sequence of times for which $F^kx=f^{t_k}x$. We have to show that
$\lim_{t\to\infty}\mbox{dist}(f^tx,\omega_f(c))=0$. Along the
subsequence $t_k$ 
this holds as
$\lim_{k\to\infty}f^{t_k}x=\lim_{k\to\infty}F^kx=c\in\omega_f(
c)$.
Consider now $t_k<t<t_{k+1}$ and suppose that $F^kx\in I_n$ (or
$F^kx\in\hat 
I_n$). As $\f{n}$ is monotone on $I_n$ (and on $\hat I_n$) and
as
$\f{n}(I_n)=\f{n}(\hat I_n)$ is an interval contained in the
union of the two 
central monotonicity intervals of $\f{n-2}$, the interval 
$V:=f^{t-t_k}(I_n)=f^{t-t_k}(\hat I_n)$ 
is contained in the union of two adjacent 
monotonicity intervals of $f^{S_{n-2}+t_{k+1}-t}$. Furthermore, 
$f^tx\in V$, and as 
$\c{n+1}\in I_n$, $V$ contains the point
$f^{S_{n+1}+t-t_k}(c)\in\omega_f(c)$.
Therefore $\mbox{dist}(f^tx,\omega_f(c))\le|V|\le
2\delta_{S_{n-2}+t_{k+1}-t}\le
2\delta_{S_{n-2}}$, where $\delta_k$ denotes the maximal length
of a 
monotonicity interval of $f^k$, and 
$\lim_{k\to\infty}\delta_k=0$ because $f$ is 
non-renormalizable. 
\qed
\par

\medskip
{\em Proof of Theorem }\ref{drift}:\quad
Let us show that we can apply Theorem \ref{random}
where we take $d=2$, $k_0=2$,
$X=(u_0,\hat u_0)$, $\A$ the partition from above
and $m$ the Lebesgue measure on $X$.
So fix $r\in \nz$ sufficiently large and consider
$A_r=I_r\cup\hat I_r$.
For $j>0$ define
\[
a_j:=|A_{r+j-3}|=2|I_{r+j-3}|=2|u_{r+j-3}-\tilde u_{r+j-2}|
\]
(observe that $f$ is symmetric),
and let
\[
a_0=\min(|c_{S_{r-2}}-u_{r-2}|,|c_{S_{r}}-u_{r}|)\ .
\]
Note that $a_1,a_2,a_3$ is the size of $A_{r-2},A_{r-1},A_r$
and that $a_0$ expresses `Koebe space'.
Now let $\mu\in\H$ be a measure of the form $\mu=F_*^n(m_{|A})$
where
$A\in\A_{n+1}$ with $F^n(A)=A_r$. Denote by $\nu$ the
normalization of $\mu$ and let for $j\ge 1$
\[
\nu_j=\nu((F_{|A_r})^{-1}(I_{r+j-3}))\ .
\]
We shall show that these numbers satisfy the {\em scaling
condition} 
provided $r$ and $\ell$ are large enough.
\par
Because of the estimates from Theorem~\ref{realb},
it follows that there exist
constants $C_1,C_2\in (0,\infty)$ such that for large
$\ell$ and large $j$,
\[
1+\frac{C_1}{\ell}\le
\frac{|u_j-c|}{|\tilde u_{j+1}-c|}\le 1+
\frac{C_2}{\ell}\ .
\]
It follows easily that for $k\ge j \ge1$,
\[
(1+\frac{C_1}{\ell})^{k-j}\le
\frac{|u_{r+j-3}-c|}{|u_{r+k-3}-c|}=
\frac{\sumai{j}}{\sumai{k}}\le
(1+\frac{C_2}{\ell})^{k-j}\
\]
provided $\ell$ and $r$ are sufficiently large.
If $j=0$ then, again because of Theorem~\ref{realb}
\[
(1+\frac{C_1}{\ell})^{k}\le
\frac{|c_{S_{r-2}}-c|}{|u_{r+k-3}-c|}=
\frac{\sumai{0}}{\sumai{k}}\le
(1+\frac{C_2}{\ell})^{k}\ .
\]
(possibly with a different constant $C_1$).
This gives condition (\ref{tfirst}).
\par
Let us now show that condition (\ref{tsecond})
is satisfied with $\jmax=2$.  
We need to estimate
\beqa\label{ger1}
\frac{a_{j+1}}{\nu_{j+1}}\frac{\nu_{k+1}}{a_{k+1}}
=
\frac{2|I_{r+j-2}|}{\nu((F_{|A_r})^{-1}(A_{r+j-2}))}
\frac{\nu((F_{|A_r})^{-1}(A_{r+k-2}))}{2|I_{r+k-2}|}
\eeqa
from below
for $j=0,1,2$, where we assume $k\ge j$.
As $A_r=I_r\cup\hat I_r$, it suffices to estimate this expression
with
$A_r$ replaced by $I_r$ and also with $A_r$ replaced by $\hat
I_r$.
Because of the symmetry of $F$, both cases can be treated in the
same way, 
and we consider without loss of generality
only the case with $I_r$. So we have to estimate
\beq\label{ger2}
\frac{|I_{r+j-2}|}{\nu((F_{|I_r})^{-1}(A_{r+j-2}))}
\frac{\nu((F_{|I_r})^{-1}(A_{r+k-2}))}{|I_{r+k-2}|}
\ .
\eeq
Let $\I$ be the partition into sets $I_i$ and $\hat I_i$, and
recall
that $F$ maps $I_r$ (and also $\hat I_r$) diffeomorphically onto
\[
\cup_{i=r}^\infty (I_i\cup \hat I_i) \cup \tilde I_{r-1}\cup
I_{r-2}\ .
\]
Denote by $\tilde I_{k}$ that one of the intervals
$I_{k}$ and $\hat I_{k}$ that is on the same side of $c$ as
$I_{r-2}$.
Then, in case that $j=0$ or $j=1$, we have
$(F_{|I_r})^{-1}(A_{r+j-2})=(F_{|I_r})^{-1}(\tilde I_{r+j-2})$
and,
as $A_{r+k-2}\supset \tilde I_{r+k-2}$, we must find 
a lower bound for
\begin{equation}\label{ger3}
\frac{|\tilde I_{r+j-2}|}{\nu((F_{|I_r})^{-1}(\tilde I_{r+j-2}))}
\frac{\nu((F_{|I_r})^{-1}(\tilde I_{r+k-2}))}{|\tilde I_{r+k-2}|}\ .
\end{equation}
Moreover, 
$I_r\subset(z_{r-1},c)\subset(c_{S_r},c)$, and
$(z_{r-1},c)$ is mapped by $F$ diffeomorphically onto
$(c_{S_{r-2}},c_{S_r})$.
As the partition
\[
\I_{n+1}=\I\vee F^{-1}(\I) \vee \dots \vee F^{-n}(\I)
\]
refines the partition $\A_{n+1}$, the set $A\in\A_{n+1}$ is a
finite union 
of intervals $H\in\I_{n+1}$
with $F^n(H)=I_r$ or $F^n(H)=\hat I_r$. Fix such an interval $H$
with $F^n(H)=I_r$.
It follows from Proposition \ref{top1} that
$F$ satisfies the following extension properties:
\begin{itemize}
\item
$F^n_{|H}$ is of the form $f^s$ (in fact, $f^s$ is a composition
of maps of the form $f^{S_i}$) and therefore $F\circ
F^n=f^{S_r+s}$;
\item
there exists an interval $T\supset H$ which is mapped by
$f^{S_r+s}$
diffeomorphically onto $(c_{S_{r-2}},c_{S_r})$.
\end{itemize}
Hence if $B$ is a subset of $I_r$, then for the measure
$\mu=F^n_*(m_{|A})$,
\[
\mu(B)=\sum_{H\in \I_{n+1}; F^n(H)=I_r}|({F^n}_{|H})^{-1}(B)|.
\]
In particular,
\beqas
&&\frac{|\tilde I_{r+j-2}|}{\nu((F_{|I_r})^{-1}(\tilde
I_{r+j-2}))}
\frac{\nu((F_{|I_r})^{-1}(\tilde I_{r+k-2}))}{|\tilde
I_{r+k-2}|}\\
&=&
\frac
{\sum_{\{H\in \I_{n+1}; F^n(H)=I_r\}}
|({F^n}_{|H})^{-1}\circ (F_{|I_r})^{-1}(\tilde I_{r+k-2})|}
{\sum_{\{H\in \I_{n+1}; F^n(H)=I_r\}}
|({F^n}_{|H})^{-1}\circ (F_{|I_r})^{-1}(\tilde I_{r+j-2})|}
\frac{|\tilde I_{r+j-2}|}{|\tilde I_{r+k-2}|}
\eeqas
which means that this last expression can be estimated from below
as the infimum over all $H\in \I_{n+1}$ with $F^n(H)=I_r$
of the expression
\[
\frac
{|({F^n}_{|H})^{-1}\circ (F_{|I_r})^{-1}(\tilde I_{r+k-2})|}
{|({F^n}_{|H})^{-1}\circ (F_{|I_r})^{-1}(\tilde I_{r+j-2})|}
\frac{|\tilde I_{r+j-2}|}{|\tilde I_{r+k-2}|}
\]
As we noted before, for each $H\in \I_{n+1}$ with $F^n(H)=I_r$,
there exists an interval $T\supset H$ such that
$F\circ F^n$ maps $T$ diffeomorphically onto
$(c_{S_{r-2}},c_{S_r})$.
Now $F$ satisfies the
assumptions of Proposition~\ref{sumlength}:
\begin{itemize}
\item the first assumption holds for obvious reasons;
\item assumptions 2), 3) and 4) of this Proposition
follow from the above extension properties
and from the bounds from Theorem~\ref{realb}
(where $\tau$ is a constant which is independent of $\ell$);
\item in assumption 5) the constant $K$ can be taken as the
intersection multiplicity $\intersect$ from Proposition~\ref{top1}.
\end{itemize}
Hence if we take an interval $T'\supset H$ such that
each component of $F\circ F^n(T'\setminus J)$
has exactly half the size of the corresponding component
of $F\circ F^n(T\setminus J)$ then
it follows from Proposition~\ref{sumlength} that
$$\sum_{m=0}^{j-1}|f^m(T')|\le K'$$
for some universal number $K'$ (provided we take $r$
sufficiently large
and therefore the set $\cup_{j\ge 0}A_{r-k_0+j}$ sufficiently small).
Hence by Proposition~\ref{crrat} the cross-ratio distortion
of $f^m|T'$ is bounded. Hence Lemma~\ref{doubc} implies that
\[
\frac
{|({F^n}_{|H})^{-1}\circ (F_{|I_r})^{-1}(\tilde I_{r+k-2})|}
{|({F^n}_{|H})^{-1}\circ (F_{|I_r})^{-1}(\tilde I_{r+j-2})|}
\frac{|\tilde I_{r+j-2}|}{|\tilde I_{r+k-2}|}
\]
is at least
\[
\O_1 \cdot
\frac{|c_{S_{r-2}}-\tilde u_{r+j-2}||c_{S_{r-2}}-\tilde
u_{r+j-1}|}
{|c_{S_{r-2}}-\tilde u_{r+k-2}||c_{S_{r-2}}-\tilde u_{r+k-1}|}
\]
if $r$ is large enough (where $\O_1$ is a universal constant).
\par
\vskip0.5cm
\hbox to \hsize{\hss\unitlength=1.3mm
\beginpic(55,110)(0,0) \let\ts\textstyle
\put(2,60){\line(1,0){62}}
\put(10,0){\line(0,1){105}}
\put(25,-2){{\it Figure 5.2.}}
\put(25,60){\squine(0, -5, -17, 0, 40, 37)}
\put(25,60){\squine(0, 6, 37, 0, -58, -53)}

\put(8.5,57.5){$c$}
\put(17.5,60){\line(0,1){30}} \put(16.5,57.5){$\tilde u_{r+1}$}
\put(45,60){\line(0,-1){50}} \put(46,57.5){$u_r$}
\put(25.5,57.5){$z_r$}
\put(57,60){\line(0,-1){53.5}} \put(51.5,57.5){$z_{r-1}$}
\put(62,60.5){\line(0,-1){1}} \put(60,57.5){$c_{S_r}$}
\put(30,60){\line(0,-1){26}} \put(31,60){\line(0,-1){30}}
\put(37,60){\line(0,-1){42}}

\put(10,10){\line(1,0){35}} \put(1.5,10){$u_{r-2}$}
\put(10,90){\line(1,0){7.5}} \put(1.5,90){$u_r$}
\put(9.5,84){\line(1,0){1}} \put(1.5,84){$\tilde u_{r+1}$}
\put(10,6.5){\line(1,0){47}} \put(1.5,6.5){$c_{S_{r-2}}$}
\put(1.5,99){$c_{S_r}$}

\put(10,34){\line(1,0){20}} \put(10,30){\line(1,0){21}}
\put(1.5,32){$I_{r+k-2}$}
\put(10,18){\line(1,0){27}}
\put(4,14){$I_{r-2}$}
\put(5,87){$I_{r}$}

\put(26,80){$f^{S_r}$}

\endpic\hss}
\vskip0.5cm
By the choice of $a_0$ this is bounded from below by
\[
\O_1 \cdot \frac{a_0}{\sumaI{k}}\cdot\frac{a_1}{\sumaI{k+1}}\
,
\]
and this is also a lower bound for (\ref{ger3}) and hence in case
$j=0$ or $1$
also for (\ref{ger2}). 
\par
If $j=2$, then both parts of 
$(F_{|I_r})^{-1}(A_{r+j-2})=
(F_{|I_r})^{-1}(I_{r+j-2})\cup (F_{|I_r})^{-1}(\hat I_{r+j-2})$
are nonempty. Admitting an additional factor of $\frac 12$ for
the lower 
bound it suffices to estimate (\ref{ger2}) with $A_{r+j-2}$ first
replaced 
by $I_{r+j-2}=I_r$ and then by $\hat I_{r+j-2}=\hat I_r$.
For $\tilde I_{r+j-2}$ (that one on the same side of $c$ as
$I_{r-2}$), 
the same estimate as above works. For $\hat{\tilde I}_{r+j-2}$ 
(that one on the other side of $c$) we estimate (\ref{ger2}) from
below
by
\begin{equation}\label{ger4}
\frac{|I_{r}|}{\nu((F_{|I_r})^{-1}(I_{r}))}
\frac{\nu((F_{|I_r})^{-1}(\hat{\tilde I}_{r+k-2}))}{|\hat{\tilde
I}_{r+k-2}|}\ ,
\end{equation}
and along the same lines as above we find the lower estimate
\[
\O_1 \cdot
\frac{|c_{S_{r}}-u_{r}||c_{S_{r}}-u_{r+1}|}
{|c_{S_{r}}-\hat{\tilde u}_{r+k-2}||c_{S_{r}}-\hat{\tilde
u}_{r+k-1}|}\ .
\]
By the choice of $a_0$ this is bounded from below by
\[
\O_1 \cdot \frac{a_0}{\sumaI{k}}\cdot\frac{a_3}{\sumaI{k+1}}\
,
\]
and this is also a lower bound for (\ref{ger2}) in case $j=2$.
This bounds
(\ref{ger1}) from below and hence proves (\ref{tsecond}).
\par
We turn to the proof of (\ref{tthird}).
\[
\frac{\nu_{k+1}}{a_{k+1}}
=
\frac{\nu_1}{a_1}
\frac{|I_{r-2}|}{\nu((F_{|A_r})^{-1}(I_{r-2}))}
\frac{\nu((F_{|A_r})^{-1}(I_{r+k-2}))}{|I_{r+k-2}|}
+\frac{\nu_3}{a_3}
\frac{|\hat I_r|} {\nu((F_{|A_r})^{-1}(\hat I_{r}))}
\frac{\nu((F_{|A_r})^{-1}(\hat I_{r+k-2}))}{|\hat I_{r+k-2}|}\
,
\]
and as above it suffices to estimate this expression with $A_r$
replaced 
by $I_r$. Using the second inequality in Lemma \ref{doubc} a
rough upperbound
for the first summand is given by
\[
\O_2 \frac{\nu_1}{a_1}\cdot
\left|\frac{c_{S_{r-2}}-c_{S_r}}{c_{S_r}-c}\right|^2
\le
\frac{\nu_1}{a_1}\cdot
\left(1+\frac{\sumaI{\infty}}{\sum_{i=3}^{\infty}a_i}\right)^2
\le
\frac{\nu_1}{a_1}\cdot
(1+\rhop^3)^2
<(27)^2\frac{\nu_1}{a_1}
\ ,
\]
and the second one similarly by
\[
\O_2 \frac{\nu_3}{a_3}\cdot
\left|\frac{c_{S_{r-2}}-c_{S_r}}{c_{S_{r-2}}-c}\right|^2
\le
\Kp\rhop^2\frac{\nu_3}{a_1}\cdot2^2
\le
16\Kp\frac{\nu_3}{a_1}\ .
\]
This yields (\ref{tthird}).

\bigskip
\end{document}